\documentclass[11pt]{article}
\usepackage{graphicx,psfrag}

\usepackage{amssymb,latexsym,amscd,amsmath,amsthm}
\usepackage[all]{xy}
\xyoption{all}
\usepackage{geometry}

\renewcommand{\square}{\hbox{${\vcenter{\hrule height.4pt
  \hbox{\vrule width.4pt height6pt \kern6pt
     \vrule width.4pt}
  \hrule height.4pt}}$}}

\renewcommand {\qed} {\hfill \nobreak \square \medbreak}

\newcommand{\be}{\begin{equation}}
      \newcommand{\ee}{\end{equation}}
      \newcommand{\ba}{\begin{eqnarray}}
       \newcommand{\ea}{\end{eqnarray}}
\newcommand{\ban}{\begin{eqnarray*}}
\newcommand{\ean}{\end{eqnarray*}}

\newcommand{\diam}{{\rm diam}}
\newcommand{\conju}{{\rm conj}}

\renewcommand{\leq}{\leqslant}
\renewcommand{\geq}{\geqslant}
\renewcommand{\le}{\leqslant}
\renewcommand{\ge}{\geqslant}

\newcommand{\Pf}{\noindent {\bf Proof:} }

\newtheorem{theorem}{Theorem}

\newtheorem{theo}{Theorem}[section]
\newtheorem{remark}[theo]{Remark}
\newtheorem{defn}[theo]{Definition}
\newtheorem{example}[theo]{Example}

\newtheorem{problem}[theo]{Open Problem}
\newtheorem{lemma}[theo]{Lemma}
\newtheorem{prop}[theo]{Proposition}
\newtheorem{coro}[theo]{Corollary}
\newtheorem{ex}[theo]{Example}

\begin{document}

\title{Conjugate Points in Length Spaces}

\author{K.\ Shankar\thanks {Partially supported by a grant from the
NSF (DMS-0513981)}\, and C.\ Sormani\thanks {Partially supported by a
PSC-CUNY Award} }

\date{}

\maketitle

\vspace{-.25cm}
\centerline{\it Dedicated to the memory of Detlef Gromoll}
\smallskip

\begin{abstract}
In this paper we extend the
concept of a conjugate point in a Riemannian manifold to geodesic
spaces.   
In particular, we introduce 
symmetric conjugate points and ultimate conjugate points
and relate these notions to prior notions developed for more
restricted classes of spaces.  We generalize the long homotopy
lemma of Klingenberg to this setting as well as the injectivity 
radius estimate
also due to Klingenberg which was used to produce closed geodesics or
conjugate points on Riemannian manifolds.  

We close with applications of these new kinds of conjugate points to
${\rm CBA}(\kappa)$ spaces: proving both known and new
theorems.  In particular we prove a Rauch comparison theorem, a
Relative Rauch Comparison Theorem, the fact that there
are no ultimate conjugate points less than $\pi$ apart
in a ${\rm CBA}(1)$ space and a few facts concerning closed
geodesics.  This paper is written to be accessible to students
and includes open problems.

\end{abstract}

\thispagestyle{empty}

\section*{Introduction}

Any student of classical Riemannian Geometry is introduced
to the notions of geodesics, exponential maps, Jacobi fields
and conjugate points.  They also learn about the cut locus and how
geodesics stop minimizing when they pass through cut or
conjugate points.  They learn how to construct
continuous families about geodesics whose endpoints are not
conjugate using the Inverse Function Theorem and how Klingenberg
used this approach to prove theorems about long homotopies
and the existence of geodesics.  Naturally when they study
spaces with sectional curvature bounds, they first learn the
Rauch comparison theorem and its implications.  Everything is
proven using geometric analysis and the fact that the
space is smooth (see for instance \cite{docarmo}).
\smallskip

On geodesic spaces things are not so
simple. There is no exponential map in general and geodesics need not
extend.  Even if they
do extend, there need not be a unique extension.  Families of
geodesics need not have Jacobi fields describing their infinitesimal
behavior.  There is not enough smoothness to apply the Inverse
Function Theorem
to prove the existence of families of geodesics. 
Rinow extended the notion of a conjugate point to
geodesic spaces with extensible geodesics and unique local
geodesics in his book \cite{Rinow}.  His notion is equivalent to what we
describe as a one sided conjugate point:  
{\em a point $q$ along a geodesic $\gamma$
emanating from $p$ is one-sided conjugate to
$p$ if there are points $q_i$ converging to $q$ and
distinct geodesics $\sigma_i, \tau_i$ joining $p$ to $q_i$ converging to
$\gamma$}.   
The same notion
was applied in the setting of ${\rm CBA}(\kappa)$
spaces by Alexander-Bishop \cite{AB-cartan} and for
exponential length spaces by the second author \cite{SorCosmos}.
\smallskip

To enable our extension of these notions to geodesic spaces,
we need to make a more precise definition of the convergence of the
geodesics (see Definition~\ref{DefGamma}). This in turn requires spaces
that are locally uniformly minimizing, i.e., spaces with neighborhoods
such that all geodesics of sufficiently small length are minimizing.
We do not require local uniqueness of geodesics as in the above
cited work.  The various definitions are given in Section 1.
\smallskip

By a theorem of Warner \cite{warner}
we see that the notion of a one-sided conjugate point
is indeed equivalent to the usual notion in the Riemannian
setting.
However, the notion of one-sided conjugate point is not the only notion
of conjugate point that is equivalent to the standard notion on a Riemannian
manifold.  While it was effective in proving the earlier results
in \cite{AB-cartan} and \cite{SorCosmos}, the
fact that $p$ does not vary restricts its applications significantly.
It cannot, for example, be used to extend the work of Klingenberg.
Jacobi fields on Riemannian manifolds are infinitesimal representations
of families of geodesics without requiring a common basepoint
and in Klingenberg's work he allows both ends of the geodesics to vary. This
leads us naturally to the notion of a \textit{symmetric conjugate point}
(Definition~\ref{defsymconj}). The definition is similar to one-sided conjugate
except now we have sequences $p_i, q_i$ converging to $p,q$ respectively
and pairs of distinct geodesics joining them converging to $\gamma$ as
before.  Nevertheless, in order to prove the long homotopy lemma 
we need an even
stronger notion of conjugate point to prove the existence of continuous
families of geodesics (Definitions \ref{defcontfam} and \ref{defcontfam2}).
This is developed in Sections 3 and 4 and leads to the notion of 
\textit{unreachable conjugate point} (Definition~\ref{defnunreach}). 
While the notion of unreachable
conjugate point is needed to prove the existence of continuous families it
is not equivalent to the Riemannian notion (see Example~\ref{schroeder}).
\smallskip


In Section~\ref{sectult} we define \textit{ultimate conjugate
points} as those that are either symmetric conjugate or unreachable
conjugate (Definition~\ref{def-ult}).
We prove this concept is an extension of the Riemannian notion
(Theorem~\ref{theo-ult-Riem}).   In fact a geodesic whose
endpoints are not ultimate conjugate has a unique family
of geodesics about it (Proposition~\ref{prop-ult3}).
It should be noted that Rinow was able to prove the existence
of a fixed basepoint continuous family away from one sided
conjugate points in his book by essentially constructing
a global exponential map \cite{Rinow}.  Alexander-Bishop were able
to prove the existence of continuous families under curvature
bounds (see Theorem~\ref{ultafterpi}).  Here we do not require the
space to have curvature bounds or an exponential map.
\smallskip

In Section~\ref{sectlonghom}
we extend the Klingenberg's Long Homotopy Lemma
to geodesic spaces.  That is we prove: {\em If a closed,
contractible, non-trivial geodesic has length less than twice the ultimate conjugate
radius, then any null homotopy $H(s,t)$ of the closed geodesic
must pass through a curve $c_{t_0}=H(\cdot,t_0)$ of length
at least twice the ultimate conjugate radius} (see Definition~\ref{ultconjrad}
and Theorem~\ref{longhomotopy}).
To prove this we first construct {\em fans} of
geodesics running along curves (see Definition~\ref{fandefn},
Lemma~\ref{fan}, Figure~\ref{fan1}).  We carefully
control the lengths of the geodesics in the fans (Lemma~\ref{lengthfan},
Figure~\ref{fan2}, Corollary~\ref{fancor})
and then apply these fans to fill in two dimensional homotopies
(Lemma~\ref{squarefan}).
\smallskip

In Section~\ref{sectcut}
we review various notions which extend the concept of the
injectivity radius of a Riemannian manifold.  Recall that on a Riemannian
manifold geodesics stop minimizing when
they are no longer the unique geodesic back to their starting point.
On a geodesic space, like a graph, a geodesic may continue
to be minimizing. To this end we offer various refined notions of
injectivity radius (and corresponding cut loci)
all of which extend the Riemannian notion: FirstInj, UniqInj, MinRad,
SymInj and UltInj. We establish various inequalites between them.
For example the unique injectivity radius (the notion in
\cite[ pg.\ 119]{BH} which is the supremum of the distance between
points where geodesics are unique)
is less than or equal to the minimal
radius (which is defined in terms of geodesics being minimal).
These radii and cut loci are related to the variety of notions defined by
Burago--Burago--Ivanov, Miller--Pak, Otsu--Shioya, Plaut and Zamfirescu. 
(c.f.\ \cite{BBI} \cite{Miller-Pak} \cite{Otsu-Shioya} \cite{Plaut-Duke92} 
\cite{Zam96}).

\smallskip

In Section~\ref{sectkling},
we prove Klingenberg's Injectivity Radius Estimate
(Theorem~\ref{kling2}): {\em  If a compact length
space, $X$, with minimal injectivity radius, ${\rm MinRad}(X)\in (0,\infty)$
then either there is a pair of ultimate conjugate points, $p,q$
with $d(p,q)\le {\rm MinRad}(X)$ or there is a closed geodesic
$\gamma:S^1\to X$ with $L(\gamma)=2{\rm MinRad}(X)$}.  It should be noted that
the minimal radius of a flat disk is infinity (Example~\ref{disk})
but that the minimal radius is finite in any compact space with
at least one extensible geodesic (Lemma~\ref{WASLEM10.5}), so
Theorem~\ref{kling2} can be used for a large class of spaces.

\smallskip


Section~\ref{sectcba}
gives a brief review of some essential points about ${\rm CBA}(\kappa)$
spaces including their definition.  
Recall that the ${\rm CBA}(\kappa)$ property is an extension to geodesic spaces
of the Riemannian notion of having sectional curvature bounded
above by $\kappa$ (Theorem~\ref{cbaRiem}).  We review work of
Gromov, Charney-Davis, Ballman-Brin and Otsu-Shioya. 
\smallskip

In Section~\ref{sectRauch}
we prove two Rauch Comparison Theorems for ${\rm CBA}(\kappa)$
spaces.  The first is a known comparison theorem
related to us by Stephanie Alexander for one-sided conjugate points and extended to
symmetric conjugate points (Theorem~\ref{Rauch}).  It states that {\em
on a ${\rm CBA}(\kappa)$ space all symmetric conjugate points are
further than $D_\kappa$ apart}.  To prove the theorem
we describe the notion of a \textit{bridge} (see Definition~\ref{def-bridge},
Figure~\ref{bridgedef}) which may be regarded as a coarse analog of
Jacobi field for length spaces. The second theorem in Section~\ref{sectRauch}
extends the Jacobi field relative comparison
theorem which says, for example, that on a Riemannian manifold
with $\textbf{sec}\le 1$, $J(t)/\sin(t)$ is nondecreasing.  Since our
Relative Rauch Comparison Theorem (Theorem~\ref{RelRauch}) is stated with
bridges and is not infinitesimal it has error terms depending on the
height of the bridges.  As the height decreases to $0$ we get the
same control on our bridges as one has on Jacobi fields
in the Riemannian setting.
\smallskip

Section~\ref{sectultafterpi} concerns
Theorem~\ref{ultafterpi}:
{\em On a ${\rm CBA}(\kappa)$ space all ultimate conjugate points are further 
than $D_\kappa$ apart}.  To prove this one must show that
any geodesic of length less than $D_{\kappa}$ in a ${\rm CBA}(\kappa)$ space
has a continuous family around it.   This is a result of Alexander-Bishop
stated within the proof of Theorem 3 in \cite{AB-cartan}.  The
proof is given in more detail in Ballman's textbook on manifolds
of nonpositive curvature (see Theorem 4.1 in \cite{ballmann95}).
A proof based on the Rauch Comparison Theorem is given here.

\smallskip

In Section~\ref{sectappl} we apply the results in Section~\ref{sectultafterpi} to the work
in the first half of the paper. First we prove Theorem~\ref{longhomCBA}: {\em
On a locally compact ${\rm CBA}(\kappa)$ space,
any null homotopy of a contractible
closed geodesic of length $<2D_\kappa$
passes through a curve of length $\ge 2D_\kappa$}.
We then prove Corollary~\ref{WASCOR10.8} of our Klingenberg Injectivity
Radius Theorem: {\em A compact ${\rm CBA}(\kappa)$ length space, $X$,
with $\kappa>0$ such that ${\rm MinRad}(X)\in (0,D_\kappa)$ has
a closed geodesic of length twice the minimal radius}.  We relate
this to the Charney-Davis proof of Gromov's systole theorem.
We close with a discussion of possible applications of
Theorem~\ref{longhomotopy} (the generalized Long Homotopy Lemma) to
${\rm CBA}(\kappa)$ spaces with higher spherical rank. While
Example~\ref{triplehemisphere} shows that one cannot hope
to conclude rigidity as in \cite{ssw}, one may still be able to prove
some interesting results for such spaces.
\smallskip

Open problems are suggested throughout the paper; see for instance
Problem~\ref{open-symnotone}, Probem~\ref{ContFam},
Problem~\ref{open-ultnotsym}, Problem~\ref{nec}, 
Problem~\ref{longhom-morse} and Problem~\ref{problem-rank1}.
Some of these are perhaps not too difficult and may be appropriate
for graduate students while others like Problem~\ref{longhom-morse}
and Problem~\ref{problem-rank1} may require significant work.

\medskip

It is a pleasure to thank Ben Schmidt for his assistance early
in the project, particularly for pointing us to 
Warner's Theorem (Theorem~\ref{warner}).
Thanks are also due to to Viktor Schroeder for pointing out
Example~\ref{schroeder} and allowing us to include it.
We are also very grateful to Stephanie Alexander for sharing with us the proof of
Theorem~\ref{Rauch}. In addition we would like to thank both
Stephanie Alexander and Dick Bishop for their graciousness in dealing
with our reproving of Theorem~\ref{ultafterpi} which we had not discovered in the literature and accidentally claimed as our own in the first version of this paper.
We also want to thank Ezra Miller and Sasha Lytchak for offering several valuable
comments and corrections after the first posting of our preprint.
Finally, we also want to thank Werner Ballmann,
Karsten Grove and Wolfgang Ziller for useful discussions.

\section{Geodesic Spaces} \label{sectgeod}

In this section we review the concept of geodesic
spaces and geodesics in such spaces.  We introduce a new kind
of convergence of geodesics in such spaces which extends the
concept of convergence of geodesics in Riemannian manifolds
and is stronger than sup norm convergence (Definition~\ref{DefGamma2}
and Lemma~\ref{DefGammaEquiv}).
We then introduce the concept of a locally uniformly minimizing
length space [Definition~\ref{localunifmin}] which includes Riemannian
manifolds and ${\rm CBA}(\kappa)$ spaces.  A review of ${\rm CBA}(\kappa)$
spaces is given in Section 9. We then prove the Geodesic Arzela Ascoli Theorem
(Theorem~\ref{arzasc}) for compact length spaces which are locally uniformly minimizing.

\begin{defn} \label{defn1}
A \textit{geodesic space}, $X$, is a
complete metric space
such that any pair of points is joined by a rectifiable curve whose
length is the distance between the points.  This curve is called a
minimizing geodesic and it is shorter than any other curve joining the two
points.
\end{defn}

These spaces are refered to as geodesic spaces in \cite{BH}
and as {strictly intrinsic geodesic spaces} in \cite{BBI}.
A complete Riemannian manifold is a geodesic space by
the Hopf-Rinow Theorem (c.f.\ \cite{docarmo}).  In fact Hopf-Rinow
proved that any complete
locally compact length space is a geodesic space (c.f. \cite{Gr-metric}).  
Note however, that other aspects of the Riemannian
Hopf-Rinow Theorem do not hold on these more general spaces.

\begin{ex}\label{linepile}
A metric space, $X$, can
 be created by taking a collection of line segments
$[0,1+1/j]$
and gluing all the left endpoints to a common point, $p$, and all the right
endpoints to a common point, $q$.  The metric on $X$ can be defined
as the infimum of the lengths of all rectifiable curves running between
the given points.  Note that $p$ and $q$ have infinitely
many geodesics running between them and the distance between them is
$1$. Although $X$ is a complete metric space, it is not a 
geodesic space  because the distance between $p$ and $q$ is not achieved.  
If we add one more segment $[0,1]$ running from $p$ to $q$ then we obtain a
geodesic space, $Y$.
\end{ex}

\begin{defn} \label{defgeod}
A \textit{geodesic} in a geodesic
space is a curve which is locally minimizing.  A \textit{closed geodesic}
is a map $\gamma:S^1\to X$ which is locally  minimizing (generalizes
smoothly closed).  All geodesics are parametrized proportional to
arclength.
\end{defn}

Geodesics are sometimes called ``local geodesics''.  Here we will
always say minimizing geodesics when the geodesic achieves the
distance between two points (following the convention
in Riemannian geometry).  Geodesics need not be unique (even locally),
they may branch and they are not necessarily extensible.  

Since there is no notion of a tangent space or a starting vector,
the notion of convergence of geodesics needs to be clarified.
In order to compare geodesic segments of different lengths, we
will reparametrize them so that they are defined on unit intervals.

\begin{defn}\label{DefGamma}
Let $\Gamma([0,1],X)\subseteq C_0([0,1],X)$ be the
space of geodesic segments, $\gamma:[0,1] \to X$, parametrized
proportional to arclength considered as a subset of the
space of continuous functions on $[0,1]$ with the
sup norm:
\be
|\gamma_1-\gamma_2|_0=\sup_{t\in [0,1]}d_X(\gamma_1(t),\gamma_2(t)).
\end{equation}
We say $\gamma_i$ converge as geodesics to $\gamma$ if they converge
in this space once they are reparametrized to be defined on $[0,1]$.
\end{defn}

Note that the sup norm is not in fact a norm but just a metric on
$C_0([0,1],X)$ and $\Gamma([0,1],X)$.

\begin{defn}\label{DefGamma2}
Define a metric $d_{\Gamma}$ on $\Gamma([0,1],X)\subseteq
C_0([0,1],X)$
to be
\be
d_{\Gamma}(\gamma_1,\gamma_2)=|\gamma_1-\gamma_2|_0 +
|L(\gamma_1)-L(\gamma_2)|
\ee
\end{defn}

See Example~\ref{cylinders} as to why this definition gives a stronger
definition of convergence than just the sup norm convergence.  Later
we also see that they give the same topology when the space is
a Riemannian manifold.

\begin{lemma} \label{geods-converge-1}
In a compact length space, if $\gamma_i$ are uniformly
minimizing in the sense that they run minimally between
points of some uniform interval of width
$2\delta$, and if they have a uniform upper bound on length, then
a subsequence converges in $C_0$ to $\gamma$ and
\begin{equation}
\lim_{i\to\infty}L(\gamma_i)=   L(\gamma).
\end{equation}
\end{lemma}

\Pf One first applies Arzela-Ascoli to show a limit $\gamma$ exists.
Since
\begin{eqnarray} \label{minineq}
\lim_{i\to\infty} d(\gamma_i(t_i-\delta),\gamma_i(t_i+\delta)) &=&
d(\gamma(t-\delta),\gamma(t+\delta)) \\
&\le & L(\gamma([t-\delta,t+\delta]))      \\
&\le & \lim_{i\to\infty} L(\gamma_i([t-\delta,t+\delta]))\\
&\le & \lim_{i\to\infty} d(\gamma_i(t_i-\delta),\gamma_i(t_i+\delta)),
\end{eqnarray}
so all are equalities, we see that $\gamma$ is a geodesic.   
To see that the lengths converge, we sum up the
lengths of the minimizing segments and take the limit again.
\qed

Without the uniform minimizing property this need not be true
as can be seen in the well known example of a cube.

\begin{ex}\label{cube}
If $X$ is a 2--dimensional cube, then squares which avoid the edges of the cube
are geodesics, but when they reach a square lying on the
edge, this last curve is not a geodesic.  The edge square is not
minimizing on any intervals about the points on the corners because
there are short cuts in the face.  Intuitively, wrap a rubber band
around four faces of a cube and then slide it to one side.  It snaps off
when you reach the edge.
\end{ex}

A study of geodesics with uniform minimizing properties appears
in \cite{SorLength}.  Here however, we do not want to restrict
our collection of geodesics and so we restrict our spaces instead.
Recall the following definition; spaces with this property have been
discussed in \cite{SorLength}.

\begin{defn} \label{unifmin}
A geodesic space is uniformly minimizing if there exists a
$\delta>0$ such that all geodesic segments of length $\leq \delta$
are minimizing.
\end{defn}

\begin{defn} \label{localunifmin}
We say a geodesic space is locally uniformly minimizing
if about every point $p$ in the space there is a neighborhood $U_p$
and a length $\epsilon_p$, such that any geodesic in the neighborhood
$U_p$ of length $\le \epsilon_p$ is minimizing.
\end{defn}

This concept is distinct from that of a locally minimizing space where 
minimizing geodesics are required to be unique and lengths are not
uniformly estimated. Note that Riemannian manifolds and ${\rm CBA}(\kappa)$
spaces are locally uniformly minimizing. See for instance
Prop 1.4 (1) and (2) on page 160 of \cite{BH} and recall that
${\rm CBA}(\kappa)$ spaces are locally ${\rm CAT}(\kappa)$.

\begin{lemma}\label{geods-converge-2}
Let $X$ be a locally uniformly minimizing geodesic space
and geodesics $\gamma_j$ converge in the sup norm to a curve $\gamma$, then
$\gamma$ is a geodesic.
\end{lemma}

\Pf
Let $U$ be a uniformly minimizing
neighborhood about $\gamma(t)$.  Take a quarter of the uniform
minimizing radius of $U$ and let $\delta<\epsilon_{\gamma(t)}/4$.
So $\gamma([t-\delta,t+\delta])$ is a curve lying in $U$.  Eventually
$\gamma_i([t-\delta,t+\delta])$ lies in $U$ so it is unique and minimizing.  
Then we apply (\ref{minineq}) 
to prove $\gamma([t-\delta,t+\delta])$ is minimizing
as well.
\qed

The cube is not locally uniformly minimizing as can be seen by taking
the point $p$ to be one of its corners. The following lemma's proof requires
the assumption of locally uniformly minimizing and uses Lemma~\ref{geods-converge-2}.
We leave the proof as an exercise.

\begin{lemma}\label{LengthCont} \label{DefGammaEquiv}
In a locally uniformly minimizing length space $X$,
if geodesics $\gamma_i$ converge to $\gamma$ in the sup norm, then
their lengths converge.
As a consequence Definitions~\ref{DefGamma} and~\ref{DefGamma2}
are equivalent on locally uniformly minimizing length spaces
including ${\rm CBA}(\kappa)$ spaces and Riemannian manifolds.
\end{lemma}



\begin{ex}\label{cylinders}
Let $X$ be the geodesic space formed by gluing an infinite
collection of cylinders of different radii $r_j$
together along a common line.  Taking $r_j=1/j$ we see this
is the isometric product of the Hawaiian earring with a line.
If we take $\gamma$ to be a geodesic of length 1 lying on the line
where all the cylinders were glued, then it is approached uniformly
by
geodesics $\gamma_j$ of length 2 lying in the jth cylinder, wrapped
extra times around.
\end{ex}

Example~\ref{cylinders} indicates the necessity of the locally
uniformly minimizing condition.  The necessity of the local
compactness can be seen in the $Y$ of Example~\ref{linepile} where
the
sequence of geodesics from $p$ to $q$ of decreasing lengths has no
limit since they are uniformly bounded away from each other
by a distance at least $1$. We can now state the Geodesic Arzela
Ascoli Theorem
which gives convergence in the sense of Definition~\ref{DefGamma2}.

\begin{theo}[Geodesic Arzela Ascoli Theorem]\label{arzasc}
In a locally uniformly minimizing locally compact geodesic space,
a sequence of geodesics with a common upper bound on their lengths
whose endpoints converge to points $p$ and $q$,
will have a subsequence which converges in $d_\Gamma$ to a
geodesic with endpoints $p$ and $q$.
\end{theo}

\Pf
Since the endpoints converge and the length is bounded above,
all the $\gamma_i$ lie in a common closed ball.  Since the space is
locally compact, this closed ball is compact.  Since the geodesics
are parametrized proportional to arclength, they are equicontinuous
and thus by Arzela--Ascoli, they have a subsequence which converges
in the sup norm to a curve.  By the local uniform compactness and
Lemmas~\ref{geods-converge-2} and~\ref{LengthCont}, the subsequence
converges in $d_\Gamma$.
\qed

We close with the statement of a well known lemma that will be useful later.

\begin{lemma}\label{tovector}
In a Riemannian manifold, $M$,
$\gamma_i$ converge to $\gamma$ in the sup norm
if and only if $\gamma_i(0)\to\gamma(0)$ and
$\gamma_i'(0)\to \gamma'(0)$ in $TM$.
\end{lemma}


In general length spaces there is no concept of vectors
or angles between geodesics emanating from a common point.
Even in ${\rm CAT}(k)$ and ${\rm CBA}(k)$ spaces, where there is a concept of
angle, the angle between two geodesics could be 0 and the
geodesics could be distinct.  In fact a tree is a length space,
and there the geodesics start together and can agree on intervals
only to diverge later.

\section{Symmetric Conjugate Points} \label{sectsym}

We can now rigorously state the definition of a one-sided conjugate point.
This notion was first introduced by Rinow \cite{Rinow}.  His definition
of conjugate point used the notion of an exponential map that he had created
using the fact that his spaces had extensible nonbranching geodesics
which were locally unique.  He then proved that his notion was
equivalent to this one (pp 414-415 \cite{Rinow}).  See the Appendix
for a translation of these pages.

\begin{defn} \label{defonesidedconj}
Given a pair of points $p$ and $q$ in a geodesic
space, $X$, we say $q$ is \textit{one-sided conjugate to} $p$ along a
geodesic $\gamma$ running between them if there exists a
sequence of $q_i$ converging to $q$ such that for each $q_i$ there
are two distinct geodesics running from $p$ to $q_i$, $\sigma_i$
and $\gamma_i$, both of which converge
to $\gamma$ as geodesics as in Definition~\ref{DefGamma2}.
By distinct we mean that there exists $t$ such that $\gamma_i(t)\neq
\sigma_i(t)$.
\end{defn}

Note that Alexander and Bishop \cite{AB-cartan} use the weaker
definition Definition~\ref{DefGamma} for the convergence
but they work in ${\rm CBA}(0)$ spaces and
the definitions are equivalent there by Lemma~\ref{DefGammaEquiv}.
In \cite{SorCosmos} the spaces were exponential length spaces, where 
geodesics have initial vectors and the vectors were required to converge. 
Zamfirescu has a stronger notion of conjugate point which requires
that the geodesics be minimizing \cite{Zam96}.
In \cite{warner}, F.\ Warner proves that Riemannian conjugate implies 
one-sided conjugate although he did not state it as such.

\begin{theo} [Warner]\label{warner}
If $M$ is a complete Riemannian manifold then $q$ is one sided conjugate 
to $p$ along $\gamma$ in the sense of Definition~\ref{defonesidedconj}
if and only if it is conjugate to $p$ in the Riemannian sense, i.e.,
there is a nontrivial Jacobi field along $\gamma$ which vanishes at
$p$ and at $q$.
\end{theo}

Warner's Theorem is easy to see in one direction just using the
Inverse Function Theorem to show that the exponential map
is locally one to one.  However, the opposite direction uses
special properties of the exponential map demonstrating that
it cannot behave like $f(x)=x^3$ which is one to one globally but has a
critical point. 

In this paper we would like to make a more natural extension of
the Riemannian definition of a conjugate point which allows for the 
variation of both endpoints.  

\begin{defn} \label{defsymconj}
Given a pair of points $p$ and $q$ in a geodesic
space, $X$, we say $q$ and $p$ are \textit{symmetrically conjugate}
along a geodesic $\gamma$ running between them if there exists a
sequence of $q_i$ converging to $q$ and $p_i$ to $p$
such that for each pair $p_i,q_i$ there are two distinct geodesics
$\sigma_i$ and $\gamma_i$ running from $p_i$ to $q_i$
both of which converge to $\gamma$ in $(\Gamma([0,1]),d_\Gamma)$.
\end{defn}

\begin{remark}\label{C_0symconj}
One could also say they are $C_0$ symmetrically conjugate if
one requires only that $\sigma_i$ and $\gamma_i$ converge to
$\gamma$ in the $C_0$ sense.
\end{remark}

Clearly if $q$ is one-sided conjugate to $p$ along $\gamma$ then
$p$ and $q$ are symmetrically conjugate along $\gamma$.  The converse
is less clear.

\begin{problem} \label{open-symnotone}
Can one construct a length space with a symmetric
conjugate point that is not a one-sided conjugate point?
\end{problem}

The above question can however be answered for Riemannian manifolds;
all
notions of conjugate are equivalent in this setting. One of these
equivalences
in the theorem below is Warner's theorem.

\begin{theo}\label{theo-sym-Riem}
Let $M$ be a complete Riemannian manifold and suppose
$p$ and $q$ are points on a geodesic $\gamma$. Then the following
statements are equivalent:
\begin{itemize}
\item[(i)]
$q$ is one-sided conjugate to $p$ along $\gamma$.

\item[(ii)]
$q$ is symmetrically conjugate to $p$ along $\gamma$.

\item[(iii)]
$q$ is conjugate to $p$ in the Riemannian sense, i.e., there
is a nontrivial Jacobi field along $\gamma$ which vanishes at $p$ and
at $q$.
\end{itemize}
\end{theo}

\Pf
If $p$ and $q$ are conjugate in the Riemannian sense then
Warner's Theorem says that $p$ is conjugate to $q$ in
the sense of Definition~\ref{defonesidedconj}, which immediately
implies they are symmetrically conjugate.

If $p$ and $q={\rm exp}_p(v)$ are not conjugate in the Riemannian sense
then the map $F:TM \to M\times M$ defined as
\begin{equation} \label{Fmap}
F(x,w)=(x, exp_x(w))
\end{equation}
has a nonsingular differential
at $(p,v)$.  Thus by the Inverse Function Theorem this
is locally invertible around $(p,q)$ back to points near
$(p,v)$.  Now suppose on the contrary that $p$ and $q$
are symmetrically conjugate, so there exists $p_i\to p$
$q_i \to q$ and $\gamma_i, \sigma_i$ distinct geodesics from $p_i$ to
$q_i$ converging to $\gamma$ as geodesics.  Lemma~\ref{tovector}
implies that
$\gamma_i'(0)$ and $\sigma_i'(0)$ are converging to $\gamma'(0)$.
Thus $F(p_i,\gamma_i'(0))=F(p_i,\sigma_i'(0))$ and it is not
locally one-to-one about $(p,v)$.
\qed

\section{Continuous Families of Geodesics}\label{sectcontfam}

Recall Definition~\ref{DefGamma} of the space of geodesics,
$\Gamma([0,1],X) \subset C_0([0,1],X)$ and $d_\Gamma$.

\begin{defn}\label{defcontfam}
We say $F$ is a continous family of geodesics about $\gamma$
if there are neighborhoods $U$ of $\gamma(0)$ and $V$ of $\gamma(1)$
such that
\begin{equation}
F:U\times V \to \Gamma([0,1],X), \quad F(x,y)=\gamma, \text{ where }
x=\gamma(0) \text{ and } y=\gamma(1),
\end{equation}
is continuous with respect to $d_{\Gamma}$.
The map $F$ is defined for any $u\in U$, $v\in V$ as $F(u,v)=\sigma$,
where
$\sigma$ is a geodesic such that $\sigma(0)=u$ and $\sigma(1)=v$.
\end{defn}

\begin{remark}
One could similarly define $C_0$ continuous families about geodesics.
\end{remark}

The following weaker definition will be essential to understanding
uniqueness of families. The crucial distinction is that here we do
note assume
continuity on the whole domain, but only at the endpoints of the
given
geodesic $\gamma$.

\begin{defn}\label{defcontfam2}
We say $F$ is a family of geodesics which is continuous at $\gamma$
if there are neighborhoods $U$ of $\gamma(0)$ and $V$ of $\gamma(1)$
such that
\begin{equation}
F:U\times V \to \Gamma([0,1],X), \quad F(x,y)=\gamma, \text{ where }
x=\gamma(0)
\text{ and } y=\gamma(1),
\end{equation}
is continuous with respect to $d_{\Gamma}$ at
$(\gamma(0),\gamma(1))$.
The map $F$ is defined for any $u\in U$, $v\in V$ as
$F(u,v)=\sigma$, where $\sigma$ is a geodesic such that
$\sigma(0)=u$ and $\sigma(1)=v$.
\end{defn}

In the above definition if there are neighborhoods $U,V$ such that
the family
is unique, i.e., for any other family of geodesics, there exist
possibly smaller
neighborhoods $U', V'$ on which the two families are the same, then
we say
that there is a unique family continuous at $\gamma$.

\begin{lemma}\label{lem-ult4}
Let $X$ be a geodesic space.  If all geodesics $\gamma$
of length $L(\gamma)<R$ have unique families of geodesics about them
which
are continuous at $(\gamma(0),\gamma(1))$, then all geodesics
$\gamma$
of length $L(\gamma)<R$ have continuous families of geodesics about
them.
\end{lemma}

\Pf
We begin with a family $F$ about a geodesic
$\gamma$ which is continuous at $(\gamma(0),\gamma(1))$.
Note that if we further restrict $U$ and $V$ we can guarantee
that all geodesics in this family have lengths $<R$ by
the definition of $d_{\Gamma}$.

We need only show $F$ is continuous on possibly smaller
neighborhoods $U$ and $V$.  If not then $F$ is not continuous at a
sequence $(p_i,q_i)$ converging to $(\gamma(0),\gamma(1))$.
Let $\gamma_i=F(p_i,q_i)$.  This means that there exists an
$\epsilon_i>0$
and $(p_{i,j},q_{i,j})$ converging to $(p_i,q_i)$ such that
\be\label{limji1}
\lim_{j\to\infty}F(p_{i,j},q_{i,j})\neq F(p_i,q_i)
\ee
where it is possible this limit does not exist.
\medskip

Since $L(\gamma_i)<R$, each $\gamma_i$ has a family $F_i$ defined
about it which is continuous at $(p_i,q_i)$ so in particular
\be
\lim_{j\to\infty}F_i(p_{i,j},q_{i,j})=\gamma_i= F(p_i,q_i).
\ee
Thus there exists a $j_i$ sufficiently large that
\be\label{limji2}
F(p_{i,j},q_{i,j})\neq F_i(p_{i,j},q_{i,j})\qquad \forall j\ge j_i.
\ee
Choosing $j_i$ possibly larger, we can also guarantee that
\be
d_{\Gamma}(F_i(p_{i,j_i},q_{i,j_i}), \gamma_i) < 1/i
\ee
and that $(p_{i,j_i},q_{i,j_i}) \to (\gamma(0),\gamma(1))$.

Thus we have two distinct geodesics $F(p_{i,j_i},q_{i,j_i})$
and $F_i(p_{i,j_i},q_{i,j_i})$ both converging to $\gamma$.
This contradicts the fact that we have a unique family that is
continuous at $\gamma$. \hfill \qed

\begin{prop} \label{RiemContFamA}
On a Riemannian manifold two points on a geodesic $\gamma$, say
$\gamma(0)$ and $\gamma(1)$ are not conjugate
along $\gamma$ if there is a continuous
family about $\gamma$ which is unique among all families
of geodesics about $\gamma$ which are continuous at $\gamma$.
\end{prop}

\Pf
If $\gamma(0)$ and $\gamma(1)$ are not conjugate
then the map $F:TM \to M\times M$ defined as
\begin{equation} \label{Fmap2}
F(x,w)=(x, \exp_x(w))
\end{equation}
is locally invertible around $(\gamma(0),\gamma(1))$ back to points
near
$(\gamma(0),\gamma'(0))$.
We then construct the
continuous family using this inverse $F^{-1}:U\times V \to TM$
followed by the exponential geodesic map: ${\rm Exp}:TM \to \Gamma([0,1])$
defined as
\be
{\rm Exp}(p,v)=\exp_p(tv).
\ee
We then apply Lemmas~\ref{tovector}
and~\ref{DefGammaEquiv} to see that ${\rm Exp} \circ F^{-1}$
provides a continuous family.  If the family is not unique
on possibly a smaller pair of neighborhoods, then
there exists $(p_i,q_i)\to (\gamma(0),\gamma(1))$  with
multiple geodesics $\gamma_i$ and $\sigma_i$ running
between them (from different families).  This implies
that $\gamma(0)$ and $\gamma(1)$ are symmetric conjugate.
By Theorem~\ref{theo-sym-Riem}, they are then Riemannian conjugate
which is a contradiction.
\qed

The following proposition is classical so we do not include
a proof.

\begin{prop} \label{RiemContFamB}
On a Riemannian manifold two points on a geodesic $\gamma$ say
$\gamma(0)$ and $\gamma(1)$ are not conjugate
along $\gamma$ if and only if there is a continuous
family about $\gamma$ which is unique among all families
of geodesics about $\gamma$ which are continuous at $\gamma$.
\end{prop}



The following theorem also holds if one uses $C_0$ symmetric
conjugate points  and $C_0$ families rather than the stronger definition
of convergence on any geodesic space $X$.  There is no need to
assume local compactness or local uniform minimality.

\begin{theo} \label{UniqueContFam}
Assume $X$ is a geodesic space.
If $x,y\in X$ are not symmetrically conjugate along $\gamma$
and there is a continuous family $F:U\times V \to \Gamma([0,1],X)$
about $\gamma$,
then there exists possibly smaller
neighborhoods $U'\subset U$ of $x$ and $V'\subset V$ of $y$ such that
$F$ restricted to $U'\times V'$
is the unique continuous family about $\gamma$ defined on $U'\times
V'$.
In fact it is the only family continuous at $(\gamma(0),\gamma(1))$
on these restricted neighborhoods.
\end{theo}

\Pf
Assume on the contrary that there exists $\delta_i\to 0$
and families
\begin{equation}
F_i: B_{\delta_i}(x)\times B_{\delta_i}(y) \to \Gamma([0,1],X)
\end{equation}
about $\gamma$ such that
\begin{equation}
\sigma_i=F_i(x_i,y_i)\neq F(x_i,y_i)=\gamma_i
\end{equation}
for some $x_i\in B_{\delta_i}(x)$ and $y_i \in B_{\delta_i}(y)$.
But then $\sigma_i(0)=x_i=\gamma_i(0)$ and
$\sigma_i(1)=y_i=\gamma_i(1)$ and by the continuity
of $F$ and $F_i$ at $(\gamma(0),\gamma(1))$
we know $\sigma_i$ and $\gamma_i$ both converge
to $\gamma$ causing a symmetric conjugate point.
\qed

We are not able to show in general that unique continuous families
exist precisely when two points are not symmetrically conjugate.

\begin{problem} \label{ContFam}
Assume $X$ is a locally compact length space which is
locally minimizing. Is it true that
two points $p,q\in X$ are not symmetrically conjugate along $\gamma$
ifand only if there are neighborhoods $U$ of $p$ and $V$ of $q$
and a unique continuous map $F:U\times V \to \Gamma([0,1],X)$
such that $F(p,q)=\gamma$ and for any $u\in U$, $v\in V$,
$F(u,v)=\sigma$, where $\sigma(0)=u$ and $\sigma(0)=v$?
\end{problem}

The main difficulty in attempting to answer this question is in
the construction of the family of geodesics.
This motivates the concept of an unreachable conjugate point
and Proposition~\ref{prop-ult2}.

\section{Unreachable Conjugate Points}\label{sectunreach}

\begin{defn}\label{defnunreach}
We say that $p$ and $q$ are unreachable conjugate points along $\gamma$
if there exist sequences $p_i, q_i \to p,q$ such that no choice of a sequence 
of geodesics $\gamma_i$ running from $p_i$ to $q_i$ converges to $\gamma$. 
Otherwise we say that a pair of points is reachable.
\end{defn}

Naturally as soon as pairs of points are reachable
there are geodesics joining nearby points which can be used to build
a family about $\gamma$ which is continuous at $\gamma$. This is the
content of the next proposition. Note that it does not require local
compactness and holds for all pairs of reachable points including symmetric 
conjugate points. The proof follows almost immediately from the definitions.

\begin{prop}\label{prop-ult2}
Assume $X$ is a geodesic space.
Suppose $\gamma(0)$ and $\gamma(1)$ are reachable
along $\gamma$.  Then there exist neighborhoods $U$ about $\gamma(0)$
and $V$ about $\gamma(1)$ and a family $F$ about $\gamma$ which
is continuous at $(\gamma(0),\gamma(1))$.
\end{prop}


\begin{ex}\label{sphereunreach}
On a standard sphere the poles are unreachable conjugate points
along any minimizing geodesic running between them.  This can be
seen by taking describing the geodesic $\gamma$
as the $0$ degree longitude
and taking $p_i$ and $q_i$ all lying on the $90$ degree longitude
where $p_i$ approach the north pole $p$ and $q_i$ approach the
south pole $q$.  Since any geodesic running between $p_i$ and $q_i$
must lie on the great circle which includes the $90$ degree longitude,
none of them is close to our given geodesic $\gamma$.

Note that these poles are also symmetric conjugate points as can be
seen by taking sequences $p_i=p$ and $q_i=q$ and $\gamma_i$ the
small positive degree longitudes and $\sigma_i$ the small negative
degree longitdues.
\end{ex}

\begin{lemma}\label{unreachtoconj}
On a Riemannian manifold, unreachable conjugate points
are Riemannian conjugate points and are therefore symmetric conjugate
points as well.
\end{lemma}


The proof follows from Proposition~\ref{RiemContFamB} and Theorem~\ref{theo-sym-Riem}.
On the other hand it is not clear what may happen in a general geodesic
space. Rinow proved that unreachable conjugate points are
one sided conjugate points on geodesic spaces with extensible
geodesics and local uniqueness by effectively constructing an
exponential map such that whenever a point is not a one sided conjugate
point (an ordinary point), the exponential map provides a
continuous family of geodesics based at the fixed point (\cite{Rinow}, pgs 414--5). In other
words, the problem below is settled for the spaces studied by Rinow which are generalizations
of Riemannian manifolds. But it is still open for locally uniformly minimizing geodesic spaces.

\begin{problem}\label{open-ultnotsym}
Does there exist a geodesic space with a geodesic $\gamma$
whose endpoints are unreachable conjugate points along $\gamma$
but are not symmetric conjugate points along $\gamma$?
See Remark~\ref{inlight}.
\end{problem}

The converse of Lemma~\ref{unreachtoconj} is not true.  In fact
already for Riemannian manifolds, symmetric conjugate does not imply unreachable
conjugate. More precisely consider the following lemma.

\begin{lemma}\label{symreach}
On a Riemannian manifold if $\gamma$ is the unique geodesic
which is minimizing from
$\gamma(0)$ to $\gamma(1)$ and these points are conjugate
along $\gamma$,
then the family of minimizing geodesics is continuous
at $\gamma(0)$ and $\gamma(1)$ so the conjugate points
are not unreachable conjugate points.
\end{lemma}


We present a specific example of such a phenomenon pointed
out to us by Viktor Schroeder.

\begin{ex}\label{schroeder}
Let $M^2$ be an ellipsoid that is not a sphere, i.e., a surface of
revolution of the curve
$\frac{x^2}{a^2} + \frac{y^2}{b^2} = 1$, where $1 < a < b$. If we
consider a generic
geodesic $\gamma$, for example any geodesic that is not parallel to
the coordinate
planes, then the conjugate locus of such a geodesic looks like a
diamond with inwardly
curving sides (c.f. \cite{docarmo}, Ch 13 Rmk 2.6). On the other
hand the cut locus is a tree.
If we consider a common vertex of the cut locus and the conjugate
locus, then we obtain
a point along $\gamma$ that satisfies the hypotheses of
Lemma~\ref{symreach}.
\end{ex}

Nevertheless the concept of unreachable conjugate points is crucial
to the existence of unique continuous families and so we define a
third
kind of conjugate point which incorporates both unreachable and
symmetric conjugate points.

\section{Ultimate Conjugate Points}\label{sectult}

\begin{defn}\label{def-ult}
Two points $p,q$ on a geodesic are said to be \textit{ultimate
conjugate points} along $\gamma$ if they are either two-sided (symmetric)
conjugate along $\gamma$ or unreachable conjugate along $\gamma$.
\end{defn}

\begin{remark}
Naturally one can define $C_0$ unreachable and $C_0$ ultimate
conjugate points.
\end{remark}

\begin{theo}\label{theo-ult-Riem}
On a Riemannian manifold points are conjugate along a given geodesic
if and only if they are ultimate conjugate along that geodesic.
\end{theo}


\begin{prop} \label{prop-ult3}
Suppose $X$ is a locally uniformly minimizing and locally compact geodesic 
space. Suppose $\gamma(0)$ and $\gamma(1)$ are not ultimate
conjugate points along $\gamma$.  Then there exists a unique
continuous family about $\gamma$ as in Definition~\ref{defcontfam}.
In fact it is the only family continuous at $\gamma$.
\end{prop}

\Pf
By Theorem~\ref{UniqueContFam} we need only show there is
a continuous family. To do this it suffices to show that the family
$F$ given in Proposition~\ref{prop-ult2} is continuous on possibly
smaller
neighborhoods $U$ and $V$.  If not then $F$ is not continuous at a
sequence $(p_i,q_i)$ converging to $(\gamma(0),\gamma(1))$.
This means that there exists an $\epsilon_i>0$ and
$(p_{i,j},q_{i,j})$ converging to $(p_i,q_i)$ such that
\be
\lim_{j\to\infty}F(p_{i,j},q_{i,j})\neq F(p_i,q_i)
\ee
where it is possible this limit does not exist.  If there
is a uniform upper bound on the lengths of a subsequence of
these geodesics then by Lemma~\ref{geods-converge-2}
a subsequence of these geodesics
converges in $C_0$ to some geodesic $\sigma_i\neq F(p_i,q_i)$ running
from $p_i$ to $q_i$.  This is where we use locally uniformly
minimizing and locally compact. We now choose $j_i$ such that
\be
d_{\Gamma}(\sigma_i, F(p_{i,j_i},q_{i,j_i}))<1/i.
\ee
This implies that $\sigma_i$ converges to
\be
\lim_{i\to\infty} F(p_{i,j_i},q_{i,j_i})=\gamma.
\ee
So does $\gamma_i=F(p_i,q_i)$.
Thus $\gamma(0)$ and $\gamma(1)$ are symmetric conjugate, and
hence ultimate conjugate, along $\gamma$, which is a contradiction.
\qed

The next proposition establishes the connection between unique
continuous
families and the non-existence of ultimate conjugate points (compare
with
Question~\ref{ContFam}). But note that its proof strongly requires
the use of
$d_{\Gamma}$ in the definition of convergence, continuity and
conjugate points.
It may not hold if one uses $C_0$ versions of the definitions. We do
not know
of any examples.

\begin{prop} \label{prop-ult4}
Suppose $X$ is a geodesic space.  If there are no ultimate
conjugate points along any geodesics of length less than $R$ then
there are unique continuous families about all geodesics $\gamma$
of length $<R$.  In fact on a sufficiently small domain
the continuous family is unique among all families continuous
at $(\gamma(0),\gamma(1))$.
\end{prop}

\Pf
By Theorem~\ref{UniqueContFam} we need only show there is
a continuous family about any geodesic $\gamma$ of length $<R$.
By Proposition~\ref{prop-ult2}, we have a family $F$ about
any such $\gamma$ which is continuous at $(\gamma(0),\gamma(1))$.
The proposition then follows from Lemma~\ref{lem-ult4}.
\qed

The following proposition is well known but we include its
proof here for completeness of exposition.  It gives us a
local existence result. Note that here we use the locally minimizing
set rather than a lack of conjugate points to obtain uniqueness.
We not not require locally uniformly minimizing.

\begin{prop} \label{propcontfama}
Suppose $X$ is a locally minimizing and locally compact
geodesic space and suppose $W$ is
a minimizing neighborhood. For any pair of points
in $W$ there is a unique minimizing geodesic joining them.
If $\gamma:[0,1]\to W$ is a geodesic, then
there exists is a continuous family $F$ extending
$\gamma$ on a sufficiently small pair of neighborhoods
about $\gamma(0)$ and $\gamma(1)$ which in fact consists of
minimizing geodesics.
\end{prop}

\Pf
Let $U,V\subset W$ be neighborhoods of $\gamma(0),\gamma(1)\in W$
respectively.  Define $F:U\times V \to \Gamma([0,1],X)$ so
that $F(u,v)$ is the unique minimizing geodesic running from
$u$ to $v$.  By construction $F(\gamma(0),\gamma(1))=\gamma$.  If
$u_i\to u$ and $v_i\to v$ then $F(u_i,v_i)=\sigma_i$ are
equicontinuous and a subsequence converges to some
$\sigma_\infty \in C([0,1],X)$ running from $u$ to $v$.
Since
\begin{equation}
d(u,v) \le L(\sigma_\infty)\le \liminf_{i\to\infty} L(\sigma_i)
\le \lim_{i\to\infty} d(u_i,v_i)=d(u,v),
\end{equation}
we know $\sigma_\infty$ is a minimizing geodesic from $u$
to $v$.  Thus it is unique and so in fact $\sigma_i$
converge to $\sigma$ as geodesics without any need for
a subsequence.  So $\lim_{i\to\infty}F(u_i,v_i)=F(u,v)$.
\qed

\section{Long Homotopy Lemma \`a la Klingenberg} \label{sectlonghom}

Now we proceed to generalize the long homotopy lemma of Klingenberg (cf.\ \cite{klingenberg})
which states that if a $c$ is a non-trivial, contractible closed geodesic
in a compact Riemannian manifold $M$ of length $\ell(c) < 2 \conju (M)$, then any
null homotopy $c_s$ of $c$ contains a curve $c_{t_0}$ of length $\ell(c_{t_0}) \geq
2 \conju M$. Given that we have a notion of conjugate point in locally uniformly
minimizing length spaces, we may now define the conjugate radius of such a
space.

\begin{defn}\label{ultconjrad}
The \textit{ultimate conjugate radius} at a point $p\in M$, denoted ${\rm
UltConj}(p)$, is the largest
value $r \in (0, \infty]$ such that along all geodesics starting at
$p$, there are no ultimate
conjugate points up to length $r$. The ultimate conjugate radius of $M$ is
\be
{\rm UltConj}(M)=\inf_{p\in M}{\rm UltConj}(p).
\ee
\end{defn}


\begin{theo}[Long Homotopy Lemma]\label{longhomotopy}
Let $M$ be a locally minimizing, geodesic space and
let $c:[0,1] \longrightarrow M$ be a non-trivial, contractible closed
geodesic of
length $\ell(c) < 2 {\rm UltConj}(M)$. Then any null homotopy $H(s,t)$ of
$c(s)$
contains a curve $c_{t_0} = H(s, t_0)$ of length $\ell(c_{t_0}) \geq
2 {\rm UltConj}(M)$.
\end{theo}

The following notion will be crucial in the proof of this theorem.

\begin{defn} \label{fandefn}
Let $M$ be a locally uniformly minimizing and locally compact
geodesic space. Let $C:[0,L]\to M$ be any continuous curve.
Then
we define the \textit{fan of $C$} as the family of geodesics,
$$
F: [0,T) \times [0,1] \rightarrow \Gamma([0,1],M), \quad F(s,t) =
\sigma_s(t), \quad
\sigma_s(0) = C(0),
$$
where $\sigma_s(t)$ is a geodesic joining the points $C(0)$ and
$C(s)$.
The map $F$ is required to be continuous with respect $d_{\Gamma}$.
\end{defn}

\begin{figure}[htbp]
\begin{center}
\includegraphics[height=4cm]{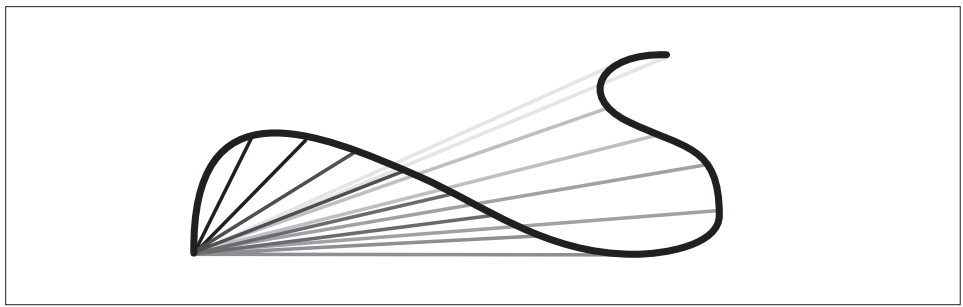}
\end{center}
\caption{A fan along an S shaped curve $C$.}
\label{fan1}
\end{figure}

The following lemma shows that for curves that are small enough, one
is able to
construct unique fans. We point out that the geodesics in the fan
need not be
minimizing, nor are they necessarily the unique geodesics joining
$C(0)$ to $C(s)$.
However, the fan is uniquely determined for $C$ as in the Lemma. See
Figure~\ref{fan1}.

\begin{lemma} \label{fan}
Let $M$ be a locally uniformly minimizing and locally compact
geodesic space.
Let $C:[0,L]\to M$ be any continuous curve parametrized by $s$. Then
there exists
a fan $F : [0,T) \times [0,1] \rightarrow \Gamma([0,1], M)$ defined
on $C$.
Furthermore, if $T$ is the maximal size of the interval on which the
fan is defined,
then either $T=L$ or $\limsup_{t\to T} L(\sigma_t)\ge {\rm UltConj}(M)$.
\end{lemma}

\Pf
We construct the fan by slowly unfolding it. Let $S$ be the
collection
of $\bar s$ such that the fan of $C$ is uniquely defined on
$[0,\bar{s}]$
and $\ell(\sigma_s)< conj_{C(0)}(M)$ for all $s \in [0,
\bar{s}]$.
Let $s_{\infty} = {\rm sup}\, S$; we need to show that either
$s_\infty = L$ or $\ell(\sigma_{s_\infty}) = {\rm UltConj}(M)$.
We know $0 \in S$, so $S$ is nonempty.  So we need only show
$S$ is closed and open as a subset of $[0,L]$.

To show $S$ is open about a given $\bar{s}\in S$
such that $\ell(\bar{s})< {\rm UltConj}(M)$,
we first extend $C$ locally below $0$ and above $L$ so that
it passes through $C(\bar{s})$. Since $C$ is a continuous curve
one may extend $C$ beyond its endpoints, for example by
backtracking on itself.

Then we just apply Proposition~\ref{prop-ult3} to the geodesic
$\sigma_{\bar{s}}$. This provides both the extension of the
definition of the
$\sigma_s$ to all $s < \bar{s}+\delta$ and the uniqueness of
such an extension.  Then we restrict ourselves back to $s\in [0,L]$,
so that we have shown $S$ is an open subset of $[0,L]$.

To prove $S$ is closed, take $s_i\in S$ converging to $s_\infty$.
Without loss of generality we may assume the $s_i$ are increasing.
So we have geodesics $\sigma_{s}$ determined by the unique
fan for all $s < s_\infty$.   Since they have a uniform upper
bound on length, they have a converging subsequence $\sigma_i$
by Theorem~\ref{arzasc}, the Geodesic Arzela--Ascoli theorem.

Let $\sigma_\infty$ be such a limit and assume
$L(\sigma_\infty)< {\rm UltConj}(M)$.
By Prop~\ref{prop-ult3}, we have a continuous family
$\tilde{\sigma}_s$ about
$\sigma_\infty$ which must include the $\sigma_i$ or
$i$ sufficiently large.  But each $\sigma_i$ also
has a unique continuous family about it (in fact this
family is $\sigma_s$), so $\sigma_s=\tilde{\sigma}_s$ near
every $s_i$. Thus $\sigma_s=\tilde{\sigma}_s$ on
$(s_\infty-\epsilon,s_\infty]$
unless there exists $t_i \to s_\infty$ where the two
families agree on one side of $t_i$ and disagree on the other side.
But this contradicts the existence of unique continuous families
about $\sigma_{t_i}=\tilde{\sigma}(t_i)$.

Once $\sigma_s=\tilde{\sigma}_s$ on $(s_\infty-\epsilon,s_\infty]$
then we see the limit $\sigma_\infty$ must have been unique and
the family is continuous on $[0,s_\infty]$ which is closed.

Since  our collection is nonempty, open and closed, it must
be all of $[0,L]$.

The fan is unique because if there were two distinct fans
then there would be another family which is continuous
at the point where they diverge from each other, and that
would imply there is a symmetric conjugate point by
Theorem~\ref{UniqueContFam}.
\qed

\begin{lemma}\label{lengthfan}
Let $X$ be a locally uniformly minimizing locally compact geodesic
space.  Let $C$ be a curve with a fan $\sigma_s$ such that
$L(\sigma_s) < {\rm UltConj}(M)$
for $s\in [0,S]$.  Then
\begin{equation}\label{lengthfaneq}
L(\sigma_s)\le L(C([0,s])).
\end{equation}
\end{lemma}

\begin{figure}[htbp]
\begin{center}
\includegraphics[height=4cm]{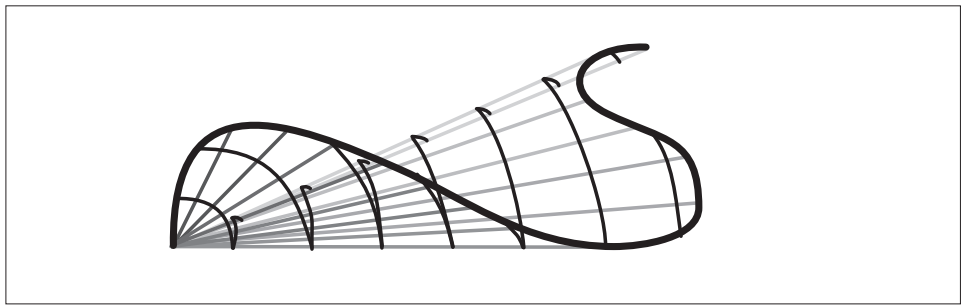}
\end{center}
\caption{A fan along an S shaped curve, $C$, with bands in black.}
\label{fan2}
\end{figure}

\Pf
We begin by defining bands $H_r \subset \bigcup_{s\in [0,S]}{\rm
Im}(\sigma_s)$
as follows; see Figure~\ref{fan2}.
\be
H_r=\{\sigma_s(r/L(\sigma_s)): s\in [0,S]\}
\ee

Note that $H_r$ is a connected set by the continuity
of the family of geodesics and that
\begin{equation}
H_r\cap {\rm Im}(C)=C(s_r)
\end{equation}
where $L(\sigma_{s_r}) = r$.

Note that $\bigcup_{s\in [0,S]} {\rm Im}(\sigma_s)$ is compact
and our space is locally uniformly minimizing, so
there exists some sufficiently small $\epsilon>0$ such that
all geodesics of length $\epsilon$ are minimizing.  In particular
all the $\sigma_s$ run minimally between bands $H_r$ and
$H_{r+\epsilon/2}$.

\begin{figure}[htbp]
\begin{center}
\includegraphics[height=4cm]{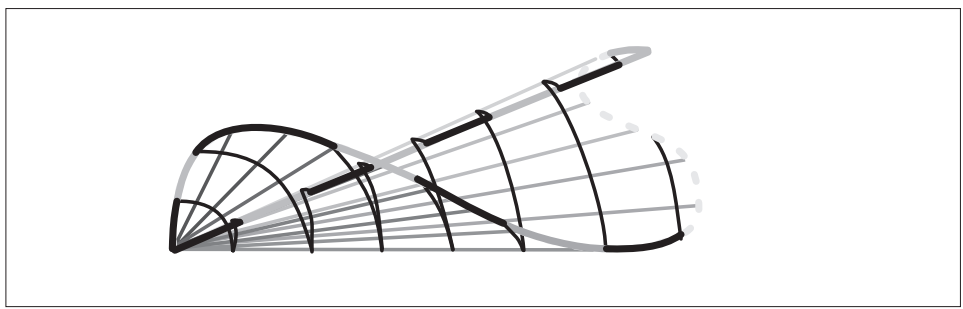}
\end{center}
\caption{The bands for each $r_i$ are marked in black, the
alternating minimizing segments along $\sigma_s$ are in black and
grey.
The corresponding segments along $C$ are marked in matching
colors (and are longer).  The total length of $C$ in fact may
include extra segments, like the dotted section, which do not
contribute to the sum.}
\label{fan3}
\end{figure}

Fixing $s \in [0,S]$ and taking a partition
$0=r_0<r_1<...<r_k=L(\sigma_s)$
such that $r_{i+1}-r_i<\epsilon/2$ we have:
\begin{eqnarray}
L(\sigma_s)&=&\sum_{i=1}^k (r_i-r_{i-1}) 
\,\,\ \le  \,\, \sum_{i=1}^k d(C(s_{r_i}),C(s_{r_{i-1}}))) \\
&\le  &\sum_{i=1}^k L(C([s_{r_i}, s_{r_{i-1}}])) 
\,\, \le \,\, L(C([0,s])).
\end{eqnarray}
Note that this last inequality follows both from the
fact that $C$ is not a geodesic and the fact that
$C$ may move back and forth in the bands.  See Figure~\ref{fan3}.
\qed

Lemma~\ref{fan} and Lemma~\ref{lengthfan} then immediately imply the
next
corollary.

\begin{coro} \label{fancor}
Let $M$ be a locally uniformly minimizing, locally compact, geodesic
space.  If $C:[0,L]\to M$ is a continuous curve of length $\le
{\rm UltConj}(M)$, then
it has a fan defined on $[0,L]$.
\end{coro}

\begin{lemma} \label{squarefan}
Let $M$ be a locally uniformly minimizing locally compact geodesic
space.  Let $H:[0,1]\times[0,1] \to M$ be continuous such that
for all $t$, $C_t(s)=H(s,t)$ implies $L(C_t)< {\rm UltConj}(M)$.  Then there
exists
a unique continuous family of geodesics
\be
\sigma:[0,1]\times[0,1]\to \Gamma([0,1],M)
\ee
such that $\forall s,t \in [0,1]$ we have
$\sigma_{s,t}(0)=H(0,0)$ and $\sigma_{s,t}(1)=H(s,t)$.
\end{lemma}

\Pf
We know by the above that for each fixed $t$ we have a unique
continuous
fan $\sigma_{t,s}$ running along $C_t$.    So we have a family
of geodesics $\sigma_{t,s}$ which is continuous with respect to the
$s$ variable.

Assume $\sigma$ is continuous on $T\times S$ where $S$ and $T$ are
intervals including $0$.  We know this is true for the trivial
intervals $S=\{0\}=T$.

We claim $S$ and $T$ are open on the right: if $t_0\in T$ and $s_0\in
S$
then there exists $\delta>0$ such that $t_0+\delta\in T$ and
$s_0+\delta \in S$.  This follows because there are
unique continuous family of geodesics about $\sigma_{t,s_0}$
and about $\sigma_{t_0,s}$ for all $s\in S$ and $t\in T$ and
so these continuous families must agree on the overlaps of all the
neighborhoods. By choosing $\delta$ small enough that
our fans land in this neighborhood, we see they must agree with
this continuous family (by the uniqueness of the fans and the
construction of the fans using the unique continuous families).

We claim $S$ and $T$ are closed on the right.  Let $s_i$ approach
${\rm sup} S=s_\infty$ and $t_i$ approach ${\rm sup} T=t_\infty$.
Then by local compactness and locally uniformly minimizing property,
and the
Geodesic Arzela Ascoli theorem, a subsequence of the
$\sigma_{t_i,s_i}$
 must converge to some limit.

Let $\sigma_\infty$ be such a limit and assume  $L(\sigma_\infty)<
\conju_{C(0)}M$.  By Proposition~\ref{prop-ult4},
we have a continuous family $\tilde{\sigma}_{s,t}$ about
$\sigma_\infty$ which must include the $\sigma_{t_i,s_i}$ for
$i$ sufficiently large.  But each $\sigma_{t_i,s_i}$ also
has a unique continuous family about it (in fact this
family is $\sigma_{s,t}$), so $\sigma_{s,t}=\tilde{\sigma}_{s,t}$
near
every $s_i$. Thus $\sigma_{s,t}=\tilde{\sigma}_{s,t}$ on
$(s_\infty-\epsilon,s_\infty]\times (t_\infty-\epsilon,t_\infty]$
unless there exists $s'_i \to s_\infty$  and $t'_i \to t_\infty$
where the two families agree on points near $(t'_i,s'_i)$ and
disagree nearby as well.
But this contradicts the existence of unique continuous families
about $\sigma_{t'_i,s'_i}=\tilde{\sigma}(t'_i,s'_i)$.

Once $\sigma_{t,s}=\tilde{\sigma}_{t,s}$ on
$(t_\infty-\epsilon,t_\infty]\times (s_\infty-\epsilon,s_\infty]$
then we see the limit $\sigma_\infty$ must have been unique and
the family is continuous on $[0,s_\infty]\times [0,t_\infty]$ which
is closed.
Since  our collection is nonempty, open and closed, it must
be all of $[0,L]$.
The family is unique because if there were two distinct families
then there would be another family which is continuous
at the point where they diverge from each other, and that would
imply there is a symmetric conjugate point by
Proposition~\ref{prop-ult4}.
\qed

\begin{figure}[htbp]
\begin{center}
\includegraphics[height=4cm]{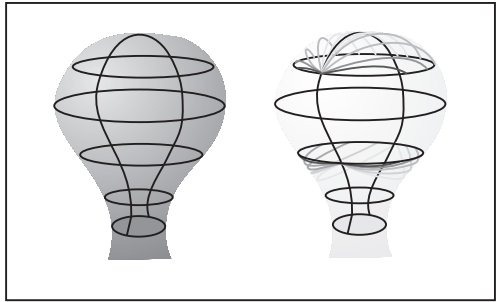}
\end{center}
\caption{The bands for each $r_i$ are marked in black, the
alternating minimizing segments along $\sigma_s$ are in black and
grey.
The corresponding segments along $C$ are marked in matching
colors (and are longer).  The total length of $C$ in fact may
include extra segments, like the dotted section, which do not
contribute to the sum.}
\label{longhom}
\end{figure}

\noindent{\bf Proof of Theorem~\ref{longhomotopy}} (Long Homotopy Lemma):
Each curve in the homotopy is denoted $c_t(s)=H(s,t)$, so
$c_0(s) =p$ (the point curve) and $c_1(s) = c(s)$ (the given
curve).
Without loss
of generality we may assume each curve $c_t$ is parametrized
proportional to arclength.
We assume, by way of contradiction, that all closed curves in the
homotopy  have length $L(c_t) < 2 {\rm UltConj}(M)$.

Let $\gamma(t) =H(0,t)$ (the starting point of each curve),
and let $\sigma(t) = H(\frac{1}{2}, t)$ (the halfway point of each
curve).
Then $\gamma(0)=\sigma(0) = p$,
$\gamma(1) = c(0)$ and $\sigma(1) = c(\frac{1}{2})$.
In Figure~\ref{longhom}: on the left we see $p$ at the top, $\gamma$
running
down the front of the surface and $\sigma$ running down the back of
the
surface, with $c$ at the neck and $c_{1/5}$, $c_{2/5}$,
$c_{3/5}$ and $c_{4/5}$  drawn as ellipses for simplicity.

By our control on the lengths of $c_t$ and
by Corollary~\ref{fancor} we have two well defined fans:
$h_{s,t}$ is the fan for $c_t(s)$ running along $s=[0,1/2]$ and
$\bar{h}_{s,t}$ which is the fan for $c_t(s)$ running backward from
$s=1$ down to $s=1/2$.  In Figure~\ref{longhom} on the right we
see two such pairs of fans, above one for $c_{1/5}$ and
below one for $c_{4/5}$.
By Lemma~\ref{squarefan},
the families $h_{s,t}$ and $\bar{h}_{s,t}$ are continuous in both
$s$ and $t$.

We will say $c_t(s)$ has a {\em closed pair of fans} when the two
fans
meet at a common geodesic,
$h_{1/2,t}(r) \equiv \bar{h}_{1/2,t}(r)$.  This happens trivially
at $t=0$, where
\be
h_{s,0}(r) \equiv\bar{h}_{s,0}(r)=H(s,0)=p.
\ee
In Figure~\ref{longhom} the pair of fans fo $c_{1/5}$
is closed
but not the pair for $c_{4/5}$.
By the continuity of $h$ and $\bar{h}$, $c_t(s)$ must
have a closed fan for all $t$.

However $c_1(s)$ is a geodesic, so when one creates a fan about
it one just gets $h_{s,1}(r)=c_1(r)$ and
$\bar{h}_{s,1}(r)=c_1(L(c)-r)$
for all values of $s$.  This is not closed.  By continuity, nearby
curves $c_t(s)$ for $t$ close to $1$ are also not
closed.  This is a contradiction.
\qed

\section{Cut Loci and Injectivity Radii}\label{sectcut}

In this section we define four kinds of cut loci and injectivity radii
all of which extend the Riemannian definition for these concepts
and yet take on different values on length spaces.  Those interested
in the ${\rm CBA}(\kappa)$ analog of the concepts we've already introduced
may skip to the Section reviewing that theory.

\begin{defn} \label{def-cut-point}
We say $q$ is a cut point of $p$ if there
are at least two distinct minimizing geodesics from
$p$ to $q$.
\end{defn}

The following definition then naturally extends the
Riemannian definition of a cut locus; see for example \cite{docarmo}.

\begin{defn} \label{def-1st-cut}
A point $q$ is in the {\em First Cut Locus} of $p$,
$q\in {\rm 1stCut}(p)$ if it is the first cut point along
a geodesic emanating from $p$.  Note that in the case
where there is a sequence of cut points $\gamma(t_i)$
of a point $p=\gamma(0)$ with $t_i$ decreasing to $t_\infty$
then $\gamma(t_\infty)$ is considered to be a first cut
point as long as there are no closer cut points along
$\gamma$, even if $\gamma(t_\infty)$ is not a cut point itself.
\end{defn}

We use the term ``first'' because we will define other
cut loci which are not based on being the first cut
point but also agree with the traditional cut locus
on Riemannian manifolds.

\begin{defn} \label{def-1st-inj}
Let the {\em First Injectivity Radius} at $p$ be defined
\be
{\rm 1stInj}(p)= d(p, {\rm 1stCut}(p))
\ee
and ${\rm 1stInj}(M)=\inf_{p\in M} {\rm 1stInj}(p)$.  If the
cut locus of a point is empty the 1st injectivity
radius of that point is infinity.
\end{defn}

The next definition is the well known notion that has appeared in various
papers. In particular, this is the notion used in the book \cite[pg.\ 119]{BH}.

\begin{defn} \label{def-uniq-inj}
The {\em Unique Injectivity Radius} of a geodesic space $X$ denoted
${\rm UniqueInj}(X)$, is the supremum over all $r\ge 0$ such that
any two points at distance at most $r$ are joined by a unique
geodesic.
\end{defn}

\begin{example}\label{disk}
On the closed flat disk the Unique Injectivity Radius is equal to $\infty$
since for any $r >0$, any two points at distance at most $r$ can be joined
by a unique geodesic (the statement is of course vacuous for $r>1$).
The first cut locus of the origin is empty for the same reason.
\end{example}

The following lemma is not too hard to show.

\begin{lemma} \label{lem-1st-uni}
On a geodesic space
\be
{\rm 1stInj}(p) = {\rm UniqueInj}(p).
\ee
\end{lemma}

The next definition of cut locus matches the definition
on Riemannian spaces given in do Carmo \cite{docarmo} and so it
is another valid extension of the concept.  It emphasizes the
minimizing properties of geodesics rather than uniqueness.

\begin{defn} \label{def-min-cut}
We say $q$ is in the {\em Minimal Cut Locus} of $p$,
$q\in {\rm MinCut}(p)$, if there is
a minimizing geodesic running from $p$ through $q$ which
is not minimizing from $p$ to any point past $q$. We are
assuming this geodesic extends past $q$.
\end{defn}

Note that the ${\rm MinCut}(p)$ is the cut locus defined
by Zamfirescu in his work on convex surfaces \cite{Zam96}.

\begin{figure}[htbp]
\begin{center}
\includegraphics[height=4cm]{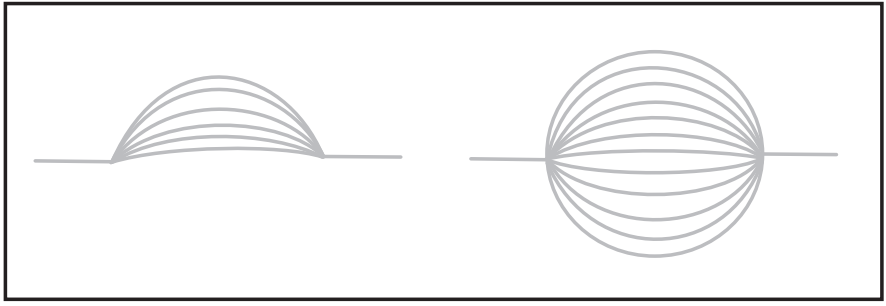}
\end{center}
\caption{The pinned sector and the pinned hemisphere.}
\label{pinsecfig}
\end{figure}

It is an essential point that this definition
has no requirement that there be a unique geodesic.
On a Riemannian manifold, geodesics stop minimizing past
the First Cut Locus because as soon as geodesics are not
unique, there are short cuts that can be taken which smooth
out the corners.  On many length spaces this is not the case:

\begin{example}\label{pinned-sector}
The pinned sector is a compact length space created by taking
an ``orange slice'' or sector between poles on a sphere and
adding line segments to each of the poles.  See
Figure~\ref{pinsecfig}.
We will denote the
poles as $p_i$ and the far ends of the
segments as $q_i$.  Then $p_2 \in UniqueCut(p_1)$
and $p_2 \notin {\rm MinCut}(p_1)$.
In contrast
the pinned closed hemisphere has $p_2 \in {\rm MinCut}(p_1)$ because the
geodesic running along the edge of the hemisphere stops minimizing at
$p_2$ and continues along the edge.  There is no such geodesic
in the pinned sector.
\end{example}

Note that while compact Riemannian
manifolds always have nonempty cut loci due to the extendability
of geodesics, this is not the case on compact length spaces:

\begin{example}
On the flat disk, the center point, $p$, has an empty {\rm Minimal
Cut Locus} because every geodesic stops existing before it
stops minimizing.  It also has an empty {\rm First Cut Locus} because
there are unique geodesics joining every point to $p$.
\end{example}

\begin{defn}\label{def-min-rad}
Let the Minimal Radius at $p$ be defined
\be
{\rm MinRad}(p)= d(p, {\rm MinCut}(p))
\ee
and ${\rm MinRad}(X)=\inf_{p\in X} {\rm MinRad}(p)$.
\end{defn}

Note that, although the minimal radius agrees with the injectivity
radius on a Riemannian manifold, we do not call it an injectivity
radius because it does not imply any uniqueness of geodesics
(like the injectivity of the exponential map).

In Example~\ref{pinned-sector}, we have 
\be
{\rm MinRad}(p)=\infty, \, {\rm 1stInj}(p)=\pi.
\ee

The following lemma follows immediately from Definition~\ref{localunifmin} of
locally uniformly minimizing and standard compactness arguments.

\begin{lemma}
A compact length space is locally uniformly minimizing if and only if its minimal 
radius is positive.
\end{lemma}

\begin{example}\label{rational-line} 
The Rationally Attached Line, is a countable collection
of real lines all attached to a common real line in the
following manner.  The first real line is attached to
the common line at corresponding integers, the second at
corresponding half integers, the third at corresponding
third integers and so on.  Note that on this line there
are infinitely many minimizing geodesics joining any
two points on the common line and that the common line
is minimizing between any pair of its values.
The unique injectivity radius and the first injectivity
radius are both zero. On the other hand, since any geodesic between
any pair of points is minimizing, the minimal radius is $\infty$.
\end{example}

In order to compare the minimal radius with the
unique injectivity radius, we recall that on a Riemannian manifold,
any point in the cut locus is either a cut point or a conjugate
point along a minimizing geodesic.  This leads to an
alternate extension of the notion of the Riemannian cut locus:

\begin{defn} \label{def-sym-cut}
The Symmetric Cut Locus of $p$, denoted ${\rm SymCut}(p)$, is the
collection of all cut points of $p$ and all symmetric
conjugate points to $p$ which are symmetric conjugate along a
minimizing geodesic.
\end{defn}

This definition is equivalent to the others on Riemannian manifolds
because the concept of a symmetric conjugate point is equivalent
to the concept of a conjugate point on Riemannian manifolds
(Theorem~\ref{theo-sym-Riem}) and on a Riemannian manifold a
symmetric conjugate point along a minimizing geodesic is a
first conjugate point. We can then also extend the definition
of injectivity radius using this notion:

\begin{defn} \label{def-sym-inj}
The Symmetric Injectivity Radius,
\be
{\rm SymInj}(p)=d(p, {\rm SymCut}(p))
\ee
and ${\rm SymInj}(X)=\inf_{p\in X} {\rm SymInj}(p)$.
\end{defn}

In Miller-Pak's work on piecewise linear spaces,
their cut locus is in fact  $Cl({\rm MinCut}(p))$.
\cite{Miller-Pak}

\begin{lemma} \label{lem-min-sym}
On a geodesic space
\be
Cl({\rm MinCut}(p)) \subset {\rm SymCut}(p)
\ee
and
\be
{\rm MinRad}(p)\ge {\rm SymInj}(p)
\ee
\end{lemma}

\Pf
Let $q_i \in {\rm MinCut}(p)$ and $q_i \to q$.  Either a subsequence of
$q_i$ are conjugate to $p$ along a minimizing geodesic $\gamma_i$
or a subsequence of $q_i$ are cut points of $p$.  In either case
there
exists pairs of geodesics $\sigma_i$ and $\bar{\sigma}_i$
running between common endpoints such that one of the two satisfies:
\be
d_{\Gamma}(\sigma_i,\gamma_i) < 1/i
\ee
and the other also satisfies this equation or has the same length
as $\gamma_i$.

Moreover, since $M$ is a uniformly locally minimizing compact space, we
can apply our Geodesic Arzela Ascoli Theorem to find a subsequence
such that $\gamma_i$ converges to $\gamma_\infty$ running
minimally between $p$ and $q$ and $\sigma_i$ converges
to $\sigma_\infty$.    If $\sigma_\infty=\gamma_\infty$,
then $p$ and $q$ are symmetric conjugate, otherwise
they are cut points.
\qed

This inequality is not an equality as can be seen by examining the
pinned sector (Example~\ref{pinned-sector}).  There the pole $p_2$
lying
opposite $p_1$ in the sector is a symmetric conjugate point for $p_1$
so $p_2 \in {\rm SymCut}(p_1)$ and
\be
{\rm SymInj}(p_1)={\rm 1stInj}(p_1)=\pi.
\ee
On the rationally attached line (Example~\ref{rational-line}),
every point is a symmetric conjugate for any other point and so
\be
{\rm SymInj}(p_1)={\rm 1stInj}(p_1)=0.
\ee
In fact the symmetric injectivity radius always agrees with the
first injectivity radius.

\begin{lemma} \label{lem-1st-sym}
On a geodesic space $X$,
\be
{\rm 1stInj}(X)= {\rm UniqueInj}(X)=
{\rm SymInj}(X) \le
{\rm MinRad}(X)
\ee
\end{lemma}

\Pf
If ${\rm UniqueInj} = r_0$, then all geodesics
of length $\le r_0$ are unique, so any sequence of pairs of
nonunique geodesics approaching a pair of symmetric conjugate points
has length at least $ r_0$, thus ${\rm SymInj}(X)\ge r_0$.

On the other hand if ${\rm SymInj}(X)=r_0$, then for any $p\in X$ we have
$d(p,{\rm SymCut}(p))\ge r_0$
and so the nearest cut point to $p$ is outside $\bar{B}_p(r_0)$,
so all geodesics of length $\le r_0$ are unique and
${\rm UniqueInj} \ge r_0$.
The rest follows from Lemma~\ref{lem-1st-uni} and
Lemma~\ref{lem-min-sym}.
\qed

Finally we introduce the ultimate cut locus and the ultimate injectivity radius.

\begin{defn} \label{def-ult-cut}
The Ultimate Cut Locus of $p$, denoted ${\rm UltCut}(p)$, is the
collection of all cut points of $p$ and all ultimate
conjugate points to $p$ which are conjugate along a
minimizing geodesic.
\end{defn}

\begin{defn} \label{def-ult-inj}
The Ultimate Injectivity Radius, denoted ${\rm UltInj}(p)$
is $d(p, {\rm UltCut}(p))$.
\end{defn}

Note that ${\rm UltInj}(p)=\min\{{\rm UltConj}(p),{\rm UniqueInj}(p)\}$
where the ultimate conjugate radius was defined in Definition~\ref{ultconjrad}.
By the definition of an ultimate conjugate point we have the following lemma.

\begin{lemma} \label{lem-sym-ult}
On a geodesic space,
\be
{\rm SymCut}(p) \subset {\rm UltCut}(p)
\ee
and
\be
{\rm UltInj}(p) \le {\rm SymInj}(p)\le {\rm MinRad}(p).
\ee
\end{lemma}

\begin{lemma} \label{lem-ultinj-fam}
If $X$ is locally uniformly minimizing and locally compact then
the ultimate injectivity radius is the supremum over all $r> 0$
such that any minimizing geodesic, $\gamma$, of length $\le r$
is unique between its endpoints and has
a continuous family about it which is unique among all families
continuous at $\gamma$ and no ultimate conjugate points along $\gamma$.
\end{lemma}

\Pf
If ${\rm UltInj}(X)=r_0$, then there are no ultimate conjugate points
along any minimizing geodesics of length less
than $r_0$ and so by Proposition~\ref{prop-ult3} there are
continuous families about these geodesics.
The uniqueness of the minimizing
geodesics of length $\le r_0$ follows from the definition
of cut and the inclusion of cut points in the ultimate cut locus.
The converse follows similarly by applying Proposition~\ref{prop-ult4}.
\qed

\begin{lemma} \label{lem-ultinj-fam-2}
If $X$ is locally uniformly minimizing and locally compact then
any geodesic, $\gamma$, of length less than the unique injectivity
radius has a continuous family about it which is unique among
all families continuous at $\gamma$.
\end{lemma}

\Pf
We will prove that the end points of $\gamma$ are not ultimate
conjugate points and apply Proposition~\ref{prop-ult4} to
complete the proof.  By Definition~\ref{def-ult} we
need only show they are neither unreachable conjugate nor
symmetric conjugate points.

First note that $\gamma(0)$ and $\gamma(1)$ are not unreachable
conjugate points. If they were, then there would be sequences of points
$p_i \rightarrow \gamma(0)$
and $q_i \rightarrow \gamma(1)$ such that any choice of sequence of
geodesics $\gamma_i$ joining $p_i$ to $q_i$ will not converge to $\gamma$.
On the other hand for $i$ large enough, $d(p_i, q_i)<{\rm UniqueInj}(X)$
since $d(\gamma(0),\gamma(1))<{\rm UniqueInj}(X)$, which
implies that there is a unique minimal geodesic, $\gamma_i$.
So $\gamma_i$ must converge in the
sup norm to $\gamma$ otherwise $\gamma(0), \gamma(1)$
would have two distinct minimizing
geodesics joining them, namely $\gamma$ and $\lim_{i \rightarrow
\infty} \gamma_i$.  Note this limit geodesic exists by the
the Geodesic Arzela-Ascoli theorem (Theorem~\ref{arzasc})
using the fact that $X$ is locally compact and locally uniformly minimizing.

We now show $\gamma(0),\gamma(1)$ are not symmetrically conjugate.
If they were then there exist sequences
$p_i \rightarrow \gamma(0)$, $q_i \rightarrow \gamma(1)$
and geodesics $\sigma_i, \lambda_i$ joining
$p_i$ and $q_i$ that converge to $\gamma$ in $d_\Gamma$.
Since $L(\gamma_i)$ and $L(\sigma_i)$ both converge to
$L(\gamma)$ eventually their lengths are
less than ${\rm UniqueInj}(X)$,
By Lemma~\ref{lem-1st-sym}, we eventually have
$$
L(\gamma_i),L(\sigma_i)< {\rm UniqueInj}(X)\le {\rm MinRad}(X)
$$
so $\gamma_i$ and $\sigma_i$ are eventually
minimizing geodesics.  However, minimizing
geodesics having length less than the Unique Injectivity
Radius are unique, so eventually $\gamma_i=\sigma_i$ which is
a contradiction.
\qed

\begin{problem}\label{nec}
It would be interesting to find an example demonstrating the
necessity of the locally uniformly minimizing condition in
the last lemma.
\end{problem}

\begin{remark}
We should point out that with any of the numerous definitions of cut loci and
injectivity radii, the injectivity radius for an arbitrary compact length space can be zero
and this is probably true generically.
\end{remark}

\section{Klingenberg's Injectivity Radius Estimate} \label{sectkling}

In the special case that $q$ is the closest point in the
minimal cut locus of $p$, Klingenberg's Lemma from
Riemannian geometry extends to locally uniformly minimizing
length spaces. 

\begin{lemma} \label{kling1}
Let $M$ be a locally uniformly minimizing, locally compact length
space. Let $p\in M$ and suppose
$q\in Cl({\rm MinCut}(p))$ such that $d(q,p)=d(p,{\rm MinCut}(p))>0$ then
either there is an ultimate conjugate point along a minimizing
geodesic from $p$ to $q$ or there is a geodesic running
from $q$ to $q$ and passing through $p$.
\end{lemma}


As you will see the proof follows
Klingenberg's original proof (c.f. \cite{docarmo}) almost exactly except that we
have Lemma~\ref{kling1}, where one shows the nearest point in the minimal cut locus
of a point either has a closer ultimate cut point or a geodesic
loop running through it.  The distinction here
is that this ultimate conjugate point need not be the closest
cut point because the geodesic running through it may
still be minimizing past it.
\medskip

\Pf
Assuming there are no ultimate conjugate points,
then by Lemma~\ref{lem-min-sym},
we know there is a pair of distinct minimizing geodesics
from $p$ to $q$: $\gamma$ and $\sigma$.  If we can show
there exists $\epsilon>0$ such that
\begin{equation}
d(\gamma(1-\epsilon),\sigma(1-\epsilon))=2\epsilon d(p,q)
\end{equation}
then we are done.

Assume, instead there exists $\epsilon_i\to 0$
such that
\begin{equation}
d(\gamma(1-\epsilon_i),\sigma(1-\epsilon_i))<2\epsilon_i d(p,q)
\end{equation}
thus the midpoint $x_i$ between these two points is closer
to $p$ than $q$ is.

Since there is a continuous family about $\gamma$ and
a continuous family about $\sigma$ we can find $\gamma_i$
converging to $\gamma$ running from $p$ to $x_i$
and we can find $\sigma_i$
converging to $\sigma$ running from $p$ to $x_i$.  Eventually
$\gamma_i$ and $\sigma_i$ must be distinct.

By the triangle inequality
\be
L(\gamma_i)< L(\gamma([0,1-\epsilon_i])+\epsilon_id(p,q) = L(\gamma)
\ee
and
\be
L(\sigma_i)< L(\sigma([0,1-\epsilon_i])+\epsilon_id(p,q) = L(\sigma).
\ee
Thus if $\gamma_i$ and $\sigma_i$ are minimizing, then we have
found a closer cut point, causing a contradiction.

By Proposition~\ref{prop-ult4} our family of geodesics
about $\gamma$ is unique among all families continuous
at $\gamma$.   Choosing $\delta>0$ sufficiently small
and $N$ sufficiently large that all the $\gamma_i$
for $i \ge N$ have $(\gamma_i(0), \gamma_i(t))$ in
the unique domain for all $t\ge t_i$ where
\be
t_i =\min \{1, (d(p,q)-\delta)/L(\gamma_i).
\ee

This implies that  $\gamma_i$ restricted
to initial segments of length $>d(p,q)-\delta$
and $\gamma$ similarly restricted must be part of this unique family.

Let us restrict $\gamma$ and $\gamma_i$
to segments of length $d(p,q)-\delta'$ for some $\delta'\in
(0,\delta)$.
The restricted $\gamma$
is minimizing and is the unique minimizing
geodesic between its end points by our definition of
${\rm Cut}(p)$.  Thus if we replace the restricted $\gamma_i$
by minimizing geodesics between their
endpoints, we'd get another family which is continuous
at $\gamma$.  By the uniqueness in Proposition~\ref{prop-ult4}
this forces the restricted $\gamma_i$ to be minimizing
between their endpoints.

That is all the $\gamma_i$ with $i\ge N$ are minimizing
on initial intervals of length $< d(p,q)$.  However,
$\gamma_i$ have length less than $d(p,q)$ so they are minimizing geodesics
without any restriction.  The same is true for $\sigma_i$ with $i$
taken sufficiently large.  Thus $x_i$ is a cut point
which leads to a contradiction.
\qed

\begin{lemma} \label{WASLEM10.5}
If $X$ is a compact length space with one
extensible geodesic $\gamma:[0,\infty)\to X$
parametrized by arclength then
\be
{\rm MinRad}(X) \le {\rm MinRad}(\gamma(0))<\infty.
\ee
\end{lemma}

\Pf
The geodesic $\gamma$ cannot be minimizing on
all intervals $[0,T]$, so eventually
it stops minimizing and we have an element
in the minimal cut locus of $\gamma(0)$.
\qed

The next theorem is applied in the final section of this paper to recover
a result of Charney-Davis on ${\rm CBA}(\kappa)$ spaces
but is much more general (Theorem~\ref{WASCOR10.8}).
Note that since the minimal injectivity radius is positive,
$X$ is locally uniformly minimizing so it satisfies the
hypothesis of Lemma~\ref{kling1}.

\begin{theo} \label{kling2}
If $X$ is a compact length space
with a positive and finite minimal injectivity radius
then either there exists $p$ and $q$
which are ultimate conjugate and
$d(p,q)\le {\rm MinRad}(X)$ or there is a closed geodesic
$\gamma:S^1\to X$ which has
length twice the
minimal injectivity radius.
\end{theo}

\Pf
Let $r_0={\rm MinRad}(X)$.  So there exists $p_i\in X$
possibly repeating such that ${\rm MinRad}(p_i)$ decreases
to $r_0$.  So $d(p_i,{\rm MinCut}(p_i))$ decreases to $r_0$.
Since $X$ is compact, there exists $q_i\in Cl({\rm MinCut}(p_i))$
such that ${\rm MinRad}(p_i)=d(p_i,q_i)$.  So by
Lemma~\ref{kling1} either $p_i$ and $q_i$ are ultimate conjugate
or there is a geodesic $\gamma_i$ running from $p_i$ to $p_i$
through $q_i$ of length $L(\gamma_i)=2d(p_i,q_i)$.

Since $X$ is compact $p_i$ has a subsequence converging to
some $p$ and $q_i$ has a subsequence converging to some $q$,
and $d(p,q)=r_0$.  Since $\gamma_i$ is minimizing on $[0,1/2]$
and on $[1/2,1]$ it has a subsequence which converges to
a piecewise geodesic $\gamma_\infty$ which is minimizing
on these intervals.

If $\gamma_\infty$ runs back and forth on the same interval
then $p$ and $q$ are symmetric conjugate and thus
ultimate conjugate.  Otherwise
$q$ is a cut point of $p$ and $q$ is a cut point of $p$.
Applying Lemma~\ref{kling1} to this pair $p$ and $q$
we see that $\gamma$ must be minimal about both
$p$ and $q$ unless they have ultimate conjugate points
on a minimizing geodesic running between them.
Therefore $\gamma:S^1 \to X$ is a geodesic.
\qed

The following example is somewhat extreme.

\begin{figure}[htbp]
\begin{center}
\includegraphics[height=6cm]{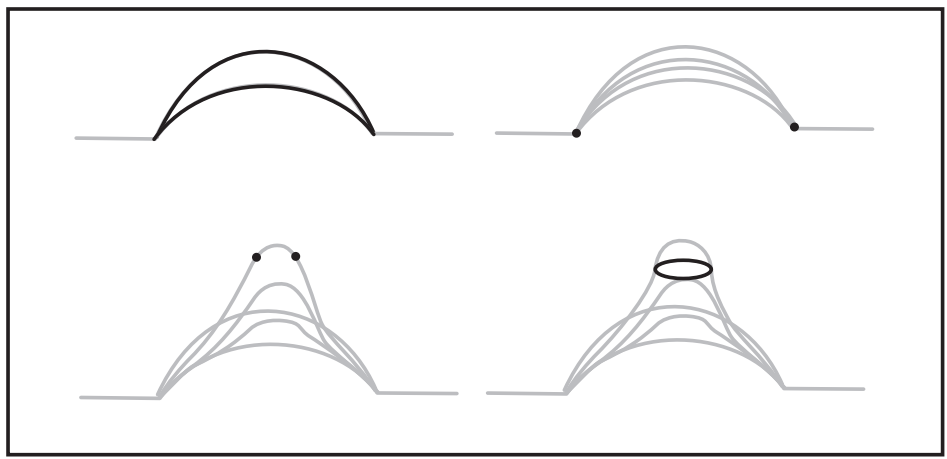}
\end{center}
\caption{The pinned sector and the pinned hemisphere.}
\label{pinsecfig2}
\end{figure}

\begin{example}
Take a circle and let $d(p_1,p_2)=r_0<\pi$.  Join $p_1$ and
$p_2$ by a second line segment of length $d$.  This space
has minimal injectivity radius and unique injectivity radius
$=r_0$.
So in this space with no ultimate conjugate points
we have the obvious geodesic loop of length $2r_0$.
See the upper left space in Figure~\ref{pinsecfig2}
where one identifies the left and right endpoints, the loop
is marked in black.
\end{example}

\begin{example}
If we fill in the loop of length $2r_0$ with a sector
of a sphere of intrinsic diameter $r_0$, then
the line segments will still be
geodesics, but the corners will provide short cuts and so the
loop will no longer be a geodesic.  In fact while the
first injectivity radius is still $r_0$ the minimal injectivity
radius
will increase to $r_0$ and there are now ultimate conjugate points on the
sector a distance $r_0$ apart.  See the upper right space
in Figure~\ref{pinsecfig2}.

One wonders if one could change the metric on the sector
so that it is still a minimizing neighborhood and somehow remove
the ultimate conjugate points so that Theorem~\ref{kling2} forces
us to see the loop which has length $\pi$.  In other words can
we introduce nonuniqueness to a circle without reducing its
minimal injectivity radius or introducing ultimate conjugate points?
In the bottom two spaces of Figure~\ref{pinsecfig2} we stretch
the sector and see that on one case we get closer ultimate
conjugate points while in the other we get a short closed loop.
\end{example}

\begin{remark}\label{inlight}
In light of the above example, one might consider using
Theorem~\ref{kling2}
to produce ultimate conjugate points in spaces where closed
geodesics are well understood.  This might help address
Open Problem~\ref{open-ultnotsym}.
\end{remark}


\section{Reviewing ${\rm CAT}(\kappa)$ versus ${\rm CBA}(\kappa)$} \label{sectcba}


In the remainder of this paper we apply the theory we have developed
to $\rm{CBA}(\kappa)$ spaces and so 
in this section we briefly recall 
the definitions of ${\rm CAT}(\kappa)$
and
${\rm CBA}(\kappa)$ spaces.   We then give a very brief review
of key results needed for this paper.  For a more in depth
approach see Bridson-Haefliger's text \cite{BH}.

Let $X$ be a metric space and $\kappa$ a fixed
real number.
To define these spaces we need the notion of a comparison
space. Let $M_\kappa^2$ denote the
2-dimensional, complete, simply connected Riemannian manifold of
constant curvature $\kappa$.  So $M_1$ is the round sphere of radius
1, $M_0$ is
the Euclidean plane and $M_{-1}$ is hyperbolic plane.  Let
$D_{\kappa}$ denote
the diameter of $M_\kappa$. Then
\be \label{defDk}
D_{\kappa} = \diam(M_\kappa) =
\begin{cases} & \frac{\pi}{\sqrt{\kappa}}, \quad \text{if $\kappa
>0$}\\
& +\infty, \quad \text{if $\kappa \leq 0$}. \end{cases}
\ee

Given a triangle constructed of three minimizing geodesics, $\Delta$,
in $X$
with perimeter less than $2D_\kappa$, we can define a comparison triangle
$\overline{\Delta}$ in $M_\kappa$ whose sides have the same length.  Given
$x$ and $y$ in $\Delta$ we can construct unique comparison points
$\bar{x}$ and $\bar{y}$ in $\bar{\Delta}$ which lie on corresponding
sides
with corresponding distances to the corresponding corners.  Since the
comparison triangle is unique up to global isometry the distance
between these comparison points does not depend on the comparison triangle.

\begin{defn}\label{defcat}
A geodesic space $X$ is a ${\rm CAT}(\kappa)$ space if for all
triangles $\Delta$ of perimeter less than $2D_\kappa$ and all points
$x,y \in \Delta$ we have
\be
d(x,y) \leq d(\bar{x}, \bar{y})
\ee
where $\bar{x}$ and $\bar{y}$ are points on the comparison triangle
$\overline{\Delta}$ in $M^2_\kappa$
\end{defn}

Note that this is a global definition; for example, a torus is not
${\rm CAT}(0)$.
To extend the concept of a Riemannian
manifold with sectional curvature bounded above one uses a localized
version of this definition:

\begin{defn} \label{defcba}
A geodesic
 space $X$ is said to be ${\rm CBA}(\kappa)$ or to have \textit{curvature
bounded
above by $\kappa$} if it is locally
a ${\rm CAT}(\kappa)$ space.  That is, for every $x \in X$, there exists
$r_x>0$ such that the metric ball
$B_x(r_x)$ with the induced metric is a ${\rm CAT}(\kappa)$ space. We
will abbreviate this to ${\rm CBA}(\kappa)$.
\end{defn}

An example of a  ${\rm CBA}(1)$ space which is not ${\rm CAT}(1)$ is
real projective
space and an example of a ${\rm CBA}(0)$ space that is not ${\rm
CAT}(0)$ is a flat torus.

\begin{theo}[Alexandrov] \label{cbaRiem}
A smooth Riemannian manifold $M$ is ${\rm CBA}(\kappa)$ if and only if
the sectional curvature of all $2$-planes in $M$ is bounded above by
$\kappa$.
\end{theo}

The condition of ${\rm CAT}(\kappa)$ is a fairly strong condition on a
space, much stronger than just having an upper curvature bound.  For example it is easy to use
the definition and degenerate comparison triangles to prove the following
well known proposition:

\begin{prop}
On a ${\rm CAT}(\kappa)$ space, $X$, all geodesics of length less than
$D_\kappa$ are minimizing and they are the unique minimizing geodesics
running between their end points:
\be
{\rm MinRad}(X)\ge D_\kappa \textrm{ and } {\rm UniqueInj}(X) \ge D_\kappa.
\ee
\end{prop}

The next proposition is then immediate implying in particular
that the results from the previous sections apply to ${\rm
CBA}(\kappa)$ spaces.

\begin{prop} \label{CBAunifmin}
A compact ${\rm CBA}(\kappa)$ space is locally uniformly minimizing as in
Definition~\ref{localunifmin}.
\end{prop}

Naturally the minimal radius of a ${\rm CBA(\kappa)}$
space may be much smaller than $D_\kappa$
just as the minimal radius of $\mathbf{RP}^2$ is only $\pi/2$ not $\pi$.
We will later show its ultimate conjugate radius is $\ge \pi$; see
Theorem~\ref{ultafterpi}.
The following theorem distinguishes between ${\rm CBA}(\kappa)$
and ${\rm CAT}(\kappa)$ spaces. Note that the notion of injectivity
radius of $X$ used by Gromov in the theorem below corresponds
to ${\rm UniqInj}(X)$ in our scheme; see Definition~\ref{def-uniq-inj}.

\begin{theo}[Gromov] \label{gromov-systole}
Let $X$ be a compact length space that is ${\rm CBA}(\kappa)$. Then $X$
fails to be ${\rm CAT}(\kappa)$
if and only if it contains a closed geodesic of length $\ell < 2D_\kappa$.
Moreover, if it contains such
a closed geodesic, then it contains a closed geodesic of length ${\rm
Sys}(X) = 2\, {\rm inj}(X)$, where the systole $Sys(X)$ is
defined to be the infimum of the lengths of all closed geodesics.
\end{theo}

In the next section of our paper we will need an angle comparison
theorem. We recall that on a ${\rm CBA}(\kappa)$ space the notion of
angle is well-defined (c.f. \cite{BH}):

\begin{defn}\label{defangle}
Let $X$ be a ${\rm CAT}(\kappa)$ space and suppose $c, c': [0,1]
\rightarrow
X$ are two geodesics issuing from the same point $p=c(0) = c'(0)$.
Then the
\textit{Alexandrov angle} between $c,c'$ at the point $p$ is defined
to be the
limit of the $\kappa$-comparison angles
$\displaystyle{\lim_{t\rightarrow 0}
\angle_p^{(\kappa)} (c(t), c'(t))}$, i.e.,
$$
\angle(c,c') = \lim_{t\rightarrow 0}\,\, 2\, {\rm arcsin}\left(
\frac{d(c(t), c'(t))}{2t}
\right)
$$
\end{defn}

Let $p,x,y$ be points in a length space $X$ such that $p\neq x, p\neq
y$.
If there are unique geodesic segments $\overline{px}, \overline{py}$,
then
we write $\angle_p (x,y)$ to denote the Alexandrov angle between
these
segments. $\angle_p^{(\kappa)} (x,y)$ denotes the angle of the
comparison triangle
in $M^2_{\kappa}$.   We then have the following
angle comparison theorem for ${\rm CAT}(\kappa)$ spaces (see for
instance
\cite{BH}, Proposition II,1.7) which may be reformulated for a suitable
${\rm CAT}(\kappa)$ neighborhood in a ${\rm CBA}(\kappa)$ space.

\begin{theo}[Alexandrov]\label{thmtopangle}
Let $X$ be a ${\rm CAT}(\kappa)$ metric space
with ${\rm MinRad}(X)\ge D_\kappa$  (i.e., there is
a unique
geodesic between any pair of points less then distance $D_\kappa$
apart).
Furthermore if $\kappa >0$, then the perimeter of each geodesic
triangle
considered is less than $2D_\kappa$.  For every geodesic triangle
$\Delta(p,q,r)$ in $X$ and for every pair of points $x \in
\overline{pq},
y\in \overline{pr}$, where $x\neq p, y\neq p$, the angles at the
vertices
corresponding to $p$ in the comparison triangles
$\overline{\Delta}(p,q,r),
\overline{\Delta}(p,x,y)$ in $M^2_\kappa$ satisfy
$$
\angle_p^{(\kappa)} (x,y) \leq \angle_p^{(\kappa)} (q,r)
$$
\end{theo}

This theorem was applied to small triangles in ${\rm CBA}(\kappa)$
spaces by Alexander-Bishop to prove the following
Cartan-Hadamard Theorem originally stated by Gromov.

\begin{theo} [Gromov]\label{gromovsystole}
If $X$ is a ${\rm CBA}(\kappa)$ length space for any $\kappa \leq 0$, then
its universal cover $\widetilde{X}$ is a ${\rm CAT}(\kappa)$ space.
\end{theo}

The original result of Gromov appears in \cite{gromov87}. A detailed
proof in the locally compact case was given by W.\ Ballmann (cf.\
\cite{ballmann90}.
Alexander--Bishop proved the theorem (cf.\ \cite{AB-cartan}) under the
additional
hypothesis that $X$ is a geodesic metric space.


Another key property of Alexandrov spaces is that they have the
nonbranching property.  That is, if two geodesics agree on an
open set then they cannot diverge later.  In particular, if there
is more than one minimal geodesic running from $p$ to $q$, then
neither geodesic can extend minimally past $q$.  As a consequence we have:

\begin{lemma}
All cut points $q$ of $p$ are in $\textrm{1stCut}(p)$ and if a minimizing
geodesic from $p$ to $q$ extends as a geodesic past $q$, then
$q\in \textrm{MinCut}(p)$ as well.
\end{lemma}

Otsu-Shioya have defined a larger cut locus, $C_p$, for Alexandrov
spaces in \cite{OtSh94} to be the collection of points $p$
which do not lie within minimal geodesics.  So $\textrm{MinCut}(p)\subset C_p$
and $\textrm{1stCut}(p) \subset C_p$.  Note that on the flat disk with $p$
at the center, $C_p$ is the boundary rather than the empty set.
Otsu-Shioya prove that the Hausdorff measure of $C_p$ is zero
\cite{Otsu-Shioya}[Prop 3.1], thus:

\begin{prop}[Otsu-Shioya]
In a ${\rm CBA}(\kappa)$ space, $\textrm{MinCut}(p)$ and
$\textrm{1stCut}(p)$ have Hausdorff measure zero.
\end{prop}

Shioya believes he can extend this to nonbranching spaces
with weaker curvature conditions like the BG scaling condition
he has been investigating recently.

\section{The Rauch Comparison Theorem}\label{sectRauch}

In this section we prove Rauch Comparison Theorems
for ${\rm CBA}(\kappa)$ spaces (Theorem~\ref{Rauch} and
Theorem~\ref{RelRauch}).  They will be applied in the
next section to prove Theorem~\ref{ultafterpi} that
ultimate conjugate points are never less than $D_\kappa$
apart and so there are continuous families about geodesics
of length $<D_\kappa$ regardless of whether they are
minimizing or not.

The following theorem and proof were outlined to us by Stephanie Alexander
for a one-sided conjugate points on ${\rm CBA}(1)$.  It is essentially
a Rauch comparison theorem.  We include its proof below for
completeness since it is necessary to motivate our extension
of the theorem, the Relative Rauch Comparison Theorem
(Theorem~\ref{RelRauch}).
Similar techniques were used by Alexander and Bishop in
their proof of the Cartan-Hadamard Theorem for ${\rm CBA}(0)$
spaces \cite{AB-cartan}.

\begin{theo}\label{Rauch}
Let $X$ be a ${\rm CBA}(\kappa)$, geodesic space. If
$\gamma(0)$ and $\gamma(1)$ are symmetric conjugate along $\gamma$,
then $L(\gamma)\ge D_\kappa$.
\end{theo}

To prove this we will build a \textit{bridge}.  In order to build the
bridge we need \textit{struts}.  We first provide precise
definitions.

\begin{defn} \label{def-strut}
A strut is a quadrilateral with a diagonal.  That is it is
a pair of minimizing geodesic segments $\sigma$ and $\gamma$ with
four sides: $S=L(\sigma)$, $T=L(\gamma)$ $A=d(\gamma(0),\sigma(0))$
$B=d(\gamma(1),\sigma(1))$ and a diagonal $D=d(\gamma(0),\sigma(1))$.
The entire strut must lie within a ${\rm CAT}(\kappa)$ neighborhood.
\end{defn}

A comparison strut can then be set up in the simply
connected two dimensional comparison space $M^2_\kappa$.  It is
built by joining two geodesic triangles with sides $SAD$ and
$TBD$.  So the comparison strut is unique when
$S+A+D$ and $T+B+D$ are less that $D_\kappa$.  The side lengths
of the comparison strut are the same but the angles all increase
by the angle version of the comparison theorem:
Definition~\ref{defangle}
and Theorem~\ref{thmtopangle}.

\begin{defn}\label{def-bridge}
A {\em bridge} on $\gamma$ is a nearby geodesic $\sigma$
and a collection of struts between them.  That is
we selected $0=s_0<s_1<....s_N=1$ and $0=t_0<t_1<...<t_N=1$
and partition $\gamma$ and $\sigma$ into segments
such that $\gamma(t_j),\gamma(t_{j+1}),\sigma(s_j),\sigma(s_{j+1})$
form a strut for $j=0...N$.

A bridge has {\em length} $\le L$ if both $\gamma$ and $\sigma$ have
length
$\le L$ and it has {\em height} $\le h$ if all struts have all their
sides
of length $\le h$.
\end{defn}

\begin{figure}[htbp]
\begin{center}
\includegraphics[height=4cm]{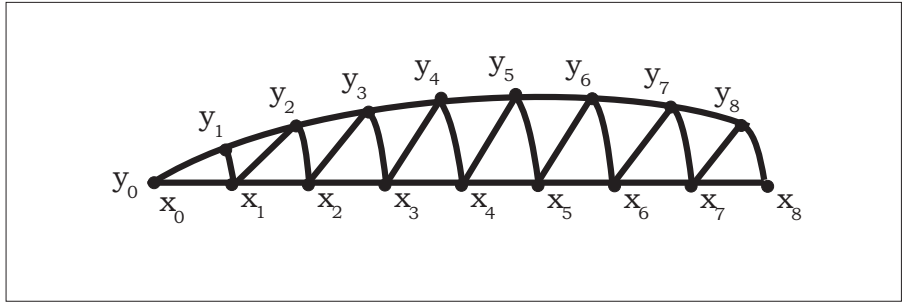}
\end{center}
\caption{A bridge with $x_j=\gamma(t_j)$ and
$y_j=\sigma(t_j)$.}
\label{bridgedef}
\end{figure}

See Figure~\ref{bridgedef}.

\begin{figure}[htbp]
\begin{center}
\includegraphics[height=6cm]{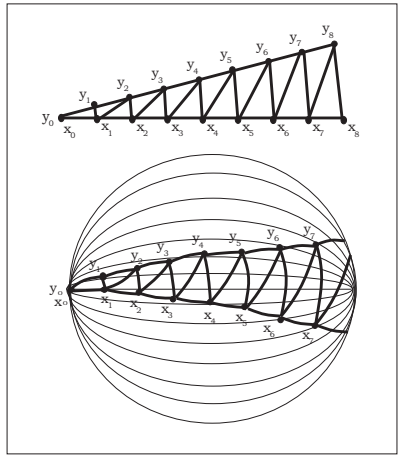}
\end{center}
\caption{Here the comparison space is a sphere.}
\label{compbridge1}
\end{figure}

One can build a comparison bridge in $M^2_\kappa$ by piecing
together the comparison struts.  Unlike the bridge in
$X$, this comparison bridge does not have geodesic beams.
The struts in fact glue together to form a bridge whose
upper and lower decks are piecewise geodesics:
$\bar{\gamma}$ and $\bar{\sigma}$.  The lengths of all the
triangles match and the interior angles in the comparison bridge
do not decrease.    See
Figure~\ref{compbridge1}.

Note it is crucial that the comparison space
is two dimensional.  In the setting where one proves a Rauch
Comparison Theorem for spaces with curvature bounded below, this leads to
some complications.  See the work of Alexander, Ghomi and Wang
\cite{AGW}.
\medskip

\noindent{\bf Proof of Theorem~\ref{Rauch}}:
Assume on the contrary that $L(\gamma)<D_\kappa$ and there
are $\gamma_i$ and $\sigma_i$
converging to $\gamma$ which share endpoints.  They
are eventually length $<T<D_\kappa$ as well by
Lemma~\ref{LengthCont} since ${\rm CBA}(\kappa)$ spaces are
locally uniformly minimizing.

Now we can cover $\gamma$ with a collection of balls such that
each ball is ${\rm CAT}(\kappa)$.  For $i$ sufficiently large both
$\gamma_i$ and $\sigma_i$ are both within the union of these
neighborhoods.  We now use $\gamma_i$ and $\sigma_i$
to build a bridge.  We do this by partitioning $\gamma_i$
and $\sigma_i$ into $N$ subsegments of equal length with
$N$ large enough that each strut formed by the corresponding
subsegments of $\gamma_i$ and $\sigma_i$ fit within a
${\rm CAT}(\kappa)$ neighborhood.

We now build a comparison bridge in $M^2_\kappa$
using the comparison struts, obtaining the piecewise geodesics:
$\bar{\gamma}_i$ and $\bar{\sigma}_i$.  Note that here we build
each triangle one by one with matching lengths starting
from $t=s=0$.  See Figure~\ref{compbridge1}.

We claim that $\bar{\gamma}_i$ and
$\bar{\sigma}_i$ are bending apart as depicted in
Figure~\ref{compbridge1}
because the interior angles are greater than $\pi$.  This can be seen
because they are built with comparison triangles.  Comparison
triangles have larger angles than the corresponding angles in
the space $X$.  Since the geodesics $\gamma$ and
$\sigma$ are straight where the triangles meet them, the
sum of the three interior angles meeting at each partition
point is $\ge \pi$, and so the sum of the comparison
angles are only larger.  Since the comparison space
is two dimensional this forces the piecewise geodesics
apart.

Thus if we were to straighten out $\bar{\gamma_i}$ and
$\bar{\sigma}_i$, creating smooth geodesics of the same length
with the same initial opening angle by bending each
geodesic segment inward, then they would meet before $1$.
That is they would meet before they are as long as $T$.  But two
distinct
geodesics cannot meet before $D_\kappa$ in $M^2_\kappa$
which gives us a contradiction.
\qed

Notice that in Theorem~\ref{Rauch} the bridge is
replacing the role of Jacobi fields in the standard differentiable
Rauch comparison theorem.  Recall that Rauch proved that a Jacobi
field $J(t)$ which has $J(0)=0$ in a space with sectional curvature less
than $\kappa=1$ then has nondecreasing
$|J(t)|/\sin(t)$.  Thus $|J(t)|\ge \sin(t)$ and it cannot hit $0$ before $\pi$.
Imitating Gromov's proof of the Relative Volume Comparison Theorem
\cite{Gr-metric}, one can say that for all $r<R<\pi$
\begin{equation} \label{diff-rel-Rauch}
\frac{|J(r)|}{\sin(r)} < \frac{|J(R)|}{\sin(R)},
\end{equation}
which we call the smooth version of the Relative Rauch Comparison
Theorem.
The inequality is opposite that of
Gromov's because the curvature is bounded above here.

We now prove the Relative Rauch Comparison Theorem
for ${\rm CBA}(\kappa)$ spaces with bridges replacing Jacobi fields.
It will be applied in the next section to prove that
there are no ultimate conjugate points before $D_\kappa$
in such spaces.

Recall Definition~\ref{def-bridge}.
First we need to insure that the bridge is short enough that
it's comparison bridge lies in a single hemisphere of the
comparison sphere.  See Figure~\ref{compbridge1}.

\begin{lemma} \label{hemisphere}
If a bridge has length $L<\pi$ and height $h\le (\pi-L)/4$
then its comparison bridge lies in a hemisphere.
\end{lemma}

\Pf
By the triangle inequality, all points on the bridge are within a
distance
$L/2 + 2h<\pi/2$ from the halfway point of one of the geodesics.
\qed

\begin{theo}[Relative Rauch Comparison] \label{RelRauch}
Let $X$ be a  ${\rm CBA}(\kappa)$
space.  Suppose $\gamma$ and $\sigma$ form a bridge in $X$
with length $L< D_\kappa$ and height $h\le (D_\kappa-L)/4$
with $\gamma(0)=\sigma(0)$.
Then for all $r\in (4h,L)$ and $R\in(r+4h, L-4h)$,
if $r_\gamma, r_\sigma \in (r-h,r+h)$ and
$R_\gamma,R_\sigma \in (R-h,R+h)$ lie
near partition points $t_i=r/L$ and $t_j=R/L$ of
the bridge,
we have:
\begin{equation} \label{bridge-rel-Rauch}
\frac{d(\gamma(r_\gamma/L),\sigma(r_\sigma/L))}
{d(\gamma(R_\gamma/L),\sigma(R_\sigma/L))}
\le
\left( 1+\frac{4h}{D_\kappa} \right)
\sup_{\bar{r} \in [r-h,r]}\frac{f_{\kappa}(\bar{r})}
{f_\kappa(R-r+\bar{r})}
+ \frac{4h+\alpha_\kappa(R-r+\bar{r},h)}{D_\kappa}
\end{equation}
where
$f_\kappa$ is the warping function of the space of
constant curvature $\kappa$
and $\alpha_\kappa(s,h)$ is an error term satisfying
\be
\lim_{h\to0}\frac {\alpha_\kappa(h,s)}{h}=0.
\ee
In particular $f_1(r)=\sin(r)$, $f_0(r)=r$,
$\alpha_1(s,h)=|s-{\rm arccos}(\cos(s)/\cos(h))|$
and $\alpha_0(s,h)=2(s-\sqrt{s^2-h^2})$.
\end{theo}

One can immediately see that as the bridge decreases in height, $h\to
0$,
our $\bar{r}$ approaches $r$ and the right hand side of
(\ref{bridge-rel-Rauch}) converges to (\ref{diff-rel-Rauch}).
Note that for $\kappa \le 0$ there is no restriction on  the
length or height of the bridge for this theorem to hold
other than requirement that the struts lie in ${\rm CAT}(\kappa)$
neighborhoods.

\begin{remark}\label{formbridge}
Note also that to apply this theorem we need only
show that $\gamma$ and $\sigma$ form a bridge with
heights $\le h/2$ in a space where all triangles of circumference
$\le 3h$ fit in ${\rm CAT}(\kappa)$ neighborhoods, and then we can always
add
extra partition points to that bridge to ensure that we pass
through $r$ and $R$.
\end{remark}

\Pf
As in the proof above
we build a comparison bridge in $M^2_\kappa$ by piecing
together the comparison struts where the struts glue together
to form a comparison bridge whose
upper and lower decks are only piecewise geodesics, $\bar{\gamma}$
and $\bar{\sigma}$, such that
\begin{equation} \label{RelRauch0}
d_{M}({\gamma}(t_j),{\sigma}(s_j))=
d_{M_\kappa}(\bar{\gamma}(t_j),\bar{\sigma}(s_j))
\end{equation}
for all $j$ including $t_j=r/L$ and $t_j=R/L$/
These piecewise geodesics bend apart because the sum
of the interior angles meeting at a given point on a deck
is greater than $\pi$ (each angle has increased in size)
exactly as depicted in Figure~\ref{compbridge1}.

This is the step where we use the ${\rm CAT}(\kappa)$.

By the triangle inequality
\begin{equation} \label{RelRauch0a}
d_{M}({\gamma}(r_\gamma/L),{\sigma}(r_\sigma/L))\le
d_{M_\kappa}(\bar{\gamma}(t_j),\bar{\sigma}(s_j))+2h
\end{equation}
and
\begin{equation} \label{RelRauch0b}
d_{M}({\gamma}(R_\gamma/L),{\sigma}(R_\sigma/L))\ge
d_{M_\kappa}(\bar{\gamma}(R/L),\bar{\sigma}(R/L))-2h
\end{equation}

Let $\bar{r}_\gamma=d_{M_\kappa}(\bar{\gamma}(r/L),\bar{\gamma}(0))$
and $\bar{r}_\sigma=d_{M_\kappa}(\bar{\sigma}(r/L),\bar{\gamma}(0))$.
Note that
\be
\bar{r}_\gamma<L(\bar{\gamma}([0,r/L]))=r \text{ and }
\bar{r}_\sigma L(\bar{\sigma}([0,r/L]))=r
\ee
because piecewise geodesics are always longer than the distance
between
their endpoints.

Since the distance is measured in $M_\kappa$ and the geodesics
$\gamma$ and $\sigma$ are bending apart towards a maximum distance
apart of $h$, we know that their length $r$ is less than the
worst path between their endpoints which runs straight and then
make a right turn and runs a distance $h$:
\be
r<\bar{r}_\gamma+h \text{ and } r<\bar{r}_\sigma+h.
\ee
 So
\be\label{wherebar}
\bar{r}_\gamma, \bar{r}_\sigma \in [r-h,r]
\ee

\begin{figure}[htbp]
\begin{center}
\includegraphics[height=6cm]{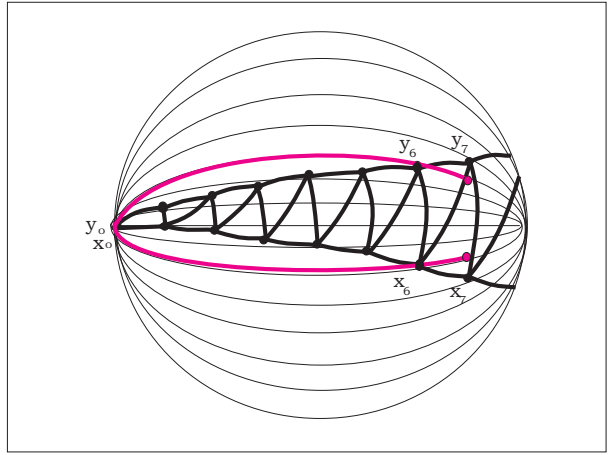}
\end{center}
\caption{Here the shaded geodesics are $\tilde{\sigma}$
and $\tilde{\gamma}$ which pass through $y_6=\bar{\sigma}(r/L)$
and $x_6=\bar{\gamma}(r/L)$.  One sees the endpoints
of $\tilde{\gamma}$ and $\tilde{\sigma}$ are closer
together than $y_7=\bar{\sigma}(R/L)$
and $x_7=\bar{\gamma}(R/L)$ since they extend the same
length from $x_6$ and $y_6$ but are bending together
relative to the comparison bridge's geodesics.}
\label{RelRauchFig3}
\end{figure}

Now we draw minimizing geodesics of length $L$, $\tilde{\gamma}$ and
$\tilde{\sigma}$ from $\bar{\gamma}(0)=\bar{\sigma}(0)$
through $\bar{\gamma}(r/L)$
and $\bar{\sigma}(r/L)$ respectively
so that $\tilde{\gamma}(\bar{r}_\gamma/L)=\bar{\gamma}(r/L)$
and $\tilde{\sigma}(\bar{r}_\sigma/L)=\bar{\sigma}(r/L)$.
 Notice
these geodesics lie outside the wedge formed by $\bar{\gamma}$
and $\bar{sigma}$ before they hit these points but then extend into
the
wedge
afterwards because the piecewise geodesics are bending apart.
See Figure~\ref{RelRauchFig3} where $\tilde{\gamma}$ and
$\tilde{\sigma}$ are in grey and $i=6<j=7$.

In particular, the points $\bar{\gamma}(R/L)$ and
$\bar{\sigma}(R/L)$ must be further apart then their
counterparts on $\tilde{\gamma}$ and $\tilde{\sigma}$
also lying a distance $R-r$ out from their common points
$\tilde{\gamma}(\bar{r}_\gamma/L)=\bar{\gamma}(r/L)$
and $\tilde{\sigma}(\bar{r}_\sigma/L)=\bar{\sigma}(r/L)$.  More
precisely:
\begin{equation} \label{RelRauch2}
d_{M_\kappa}( \bar{\gamma}((R-r+\bar{r}_\gamma)/L),
\bar{\gamma}((R-r+\bar{r}_\sigma)/L)) )\ge
d_{M_\kappa}(\tilde{\gamma}(R/L),\tilde{\sigma}(R/L) ).
\end{equation}
Here we have used the restriction on the height and length
of the bridge when $\kappa>0$ so that we avoid passing through
a pole or having a short cut running between
$\bar{\gamma}(R/L)$ $\bar{\gamma}(R/L)$ on the opposite side of the
sphere.

Applying (\ref{wherebar}) and selecting any
\be\label{wherebar2}
\bar{r} \in [r-h,r]
\ee
the triangle inequality combined with (\ref{RelRauch2})
gives us
\begin{equation} \label{RelRauch2a}
d_{M_\kappa}(\bar{\gamma}(R/L),\bar{\sigma}(R/L))+2h \ge
d_{M_\kappa}(\tilde{\gamma}((R-r+\bar{r})/L),\tilde{\gamma}((R-r+\bar{r})/L)
)
\end{equation}
and the triangle inequality combined with the choice of
$\tilde{\gamma}$
and $\tilde{\sigma}$ gives us
\begin{equation}
d_{M_\kappa}(\bar{\gamma}((r/L),\bar{\gamma}(r/L))-2h \le
d_{M_\kappa}(\tilde{\gamma}(\bar{r}_\gamma/L),
\tilde{\sigma}(\bar{r}_\sigma/L)).
\end{equation}
Thus
\be
\frac
{d_{M_\kappa}(\bar{\gamma}(R/L),\bar{\sigma}(R/L))+2h}{d_{M_\kappa}
(\bar{\gamma}(r/L),
\bar{\gamma}(r/L))-2h}
\ge
\frac
{d_{M_\kappa}(\tilde{\gamma}((R-r+\bar{r})/L),\tilde{\gamma}((R-r+\bar{r})/L))}
{d_{M_\kappa}(\tilde{\gamma}(\bar{r}/L),\tilde{\sigma}(\bar{r}/L)}.
\ee
Combining this with (\ref{RelRauch0a}) and (\ref{RelRauch0b})
we have
\be \label{RelRauch1a}
\frac
{d_{M}(\gamma(R_{\gamma}/L),\sigma(R_{\sigma}/L))+4h}
{d_{M}( \gamma(r_{\gamma}/L),\sigma(r_{\sigma}/L) )-4h}
\ge
\frac
{d_{M_\kappa}(\tilde{\gamma}((R-r+\bar{r})/L),\tilde{\sigma}((R-r+\bar{r})/L))}
{d_{M_\kappa}(\tilde{\gamma}(\bar{r}/L),\tilde{\sigma}(\bar{r}/L)) }.
\ee

Note that both $\tilde{\gamma}$ and $\tilde{\sigma}$ are
minimizing geodesics since $L<D_\kappa$
and they lie on $M_\kappa$. Thus there distances can be
computed using the warping function $f_\kappa$:

If we take $d'_{M_\kappa}$ to be the distance in $M_\kappa$
measured by taking arc paths in spheres about
$\tilde{\gamma}(0)=\tilde{\sigma}(0)$ then:
\begin{equation} \label{RelRauch3'}
d'_{M_\kappa}(\tilde{\gamma}((R-r+\bar{r})/L),\tilde{\sigma}((R-r+\bar{r})/L
)
)
= d'_{M_\kappa}(\tilde{\gamma}(\bar{r}/L),\tilde{\sigma}(\bar{r}/L))
H
\end{equation}
where
\be \label{supH}
H=\frac{f_\kappa(R-r+\bar{r})}{f_\kappa(\bar{r})}\le
\sup_{\bar{r} \in [r-h,r]}\frac{f_{\kappa}(\bar{r})}
{f_\kappa(R-r+\bar{r})}
\ee
Actual distances $d_{M_\kappa}$ are somewhat shorter but can be
estimated by this $d_{M_\kappa}'$:
\be
d_{M_\kappa}(\tilde{\gamma}(s),\tilde{\sigma}(s))
\le d_{M_\kappa}'(\tilde{\gamma}(s),\tilde{\sigma}(s))
\le d_{M_\kappa}(\tilde{\gamma}(s),\tilde{\sigma}(s)) +
\alpha_\kappa(h,s)
\ee
where $\alpha_\kappa(h,s)= 2(s-s')$ where $s-s'$ is the maximal
distance
between the arc joining $\tilde{\gamma}(s)$ and $\tilde{\sigma}(s))$
and the minimizing geodesic between these curves.  On Euclidean
space we have $s'=\sqrt{s^2-h^2}$ and on $S^2$ we have
$\cos(s) = \cos(h) \cos(s')$.
In general we can say that
\be
\lim_{h\to 0} \frac{\alpha_\kappa(s,h)}{h}=0.
\ee
Thus (\ref{RelRauch3'}) implies
\begin{equation} \label{RelRauch3b}
d_{M_\kappa}(\tilde{\gamma}((R-r+\bar{r})/L),\tilde{\sigma}((R-r+\bar{r})/L)
)
\ge d_{M_\kappa}(\tilde{\gamma}(\bar{r}/L),\tilde{\sigma}(\bar{r}/L))
H
-\alpha_\kappa(R-r+\bar{r},h)
\end{equation}
Putting this together with (\ref{RelRauch1a})
we obtain
\begin{equation} \label{RelRauch1b}
\frac{d(\gamma(r_{\gamma}/L),\sigma(r_{\sigma}/L))-4h}
{d(\gamma(R_{\gamma}/L),\sigma(R_{\sigma}/L))+4h
+ \alpha_\kappa(R-r+\bar{r},h)}
\le
\sup_{\bar{r} \in [r-h,r]}\frac{f_{\kappa}(\bar{r})}
{f_\kappa(R-r+\bar{r})}
\end{equation}
Since $d(\gamma(R_{\gamma}/L),\sigma(R_{\sigma}/L))\le D_\kappa$ and
\be
(x-4h)/(y+4h+\alpha) \le H
\text{ implies } x \le Hy +4hH +\alpha H +4h
\ee
implies
\be
x/y \le H +(4hH+\alpha h+4h)/y \le (1+4h/D_\kappa)H +
(4h+\alpha)/D_\kappa
\ee
we have (\ref{bridge-rel-Rauch}).
\qed

\section{Ultimate Conjugate Points after $\pi$} \label{sectultafterpi}

In this section we give
 a proof of the following theorem as an application of our Relative
Rauch Comparison Theorem [Theorem~\ref{RelRauch}].  We will apply it combined with our
Long Homotopy Lemma [Theorem~\ref{longhomotopy}] to prove Lemma~\ref{longhomCBA}.
We then combine it with Klingenberg's Injectivity Radius Estimate [Theorem~\ref{kling2}]
to prove a result of Charney-Davis used
in their proof of Gromov's Systole Theorem [Corollary~\ref{WASCOR10.8}].

\begin{theo}[Alexander--Bishop] \label{ultafterpi}
If $M$ is a ${\rm CBA}(\kappa)$ space and
$\gamma:[0,1]\to M$ is a geodesic of
length $< D_\kappa$, then
$\gamma$ has a unique continuous family of
geodesics about it.
\end{theo}

One may wish to keep the following example in mind when considering this theorem
in contrast to the usual smooth projective spaces:

\begin{ex}[Tetrahedral Bi-sphere]\label{tetrasphere}
Suppose two copies of the round sphere $\mathbf{S}^2$ are
glued
to each other at $4$ points that form the vertices of a regular
tetrahedron
on each sphere. The resulting space is ${\rm CBA}(1)$, but it does
not
have unique minimizing geodesics between all pairs of points of
distance
less than $\pi$.
\end{ex}

Note that by Theorem~\ref{Rauch}, we already
know there are no symmetric conjugate points along $\gamma$
of length less than $D_\kappa$.
So if we show every geodesic of length less than $\pi$ has
a continuous family about it (and thus no unreachable conjugate
points) then there are no ultimate conjugate points and,
by Proposition~\ref{prop-ult4}, we know that there is a unique
continuous family about $\gamma$.  
In fact one can see that Alexander-Bishop
proved exactly what we require in \cite{AB-cartan} and it is applied in
\cite{AB-diffgeo}.  A precise statement of their result and proof can
be also be found in Ballman's book on manifolds with nonpositive curvature
\cite{ballmann95}[Ch 1 Thm 4.1]),
as pointed out to us by Lytchak.

For completeness of exposition we provide
a proof based on the Relative Rauch 
Comparison Theorem [Theorem~\ref{RelRauch}].

In our proof, we
assume there is a continuous family about $\gamma$ restricted
to some interval $[0,T]$ and extend that family out
to $[0,T+T_1]$ for $T_1>0$ but close enough that
$\gamma([T-T_1,T+T_1])$ is in a minimizing neighborhood
so that we can apply Lemma~\ref{lem-ultinj-fam-2} and
Proposition~\ref{CBAunifmin}].
To do this we glue geodesics in the continuous family
to geodesics in the minimizing neighborhood.  However it
is not clear how one can glue geodesics to form geodesics rather than
just piecewise geodesics.  In fact, this is impossible when
$T=D_\kappa$.

\begin{figure}[htbp]
\begin{center}
\includegraphics[height=6cm]{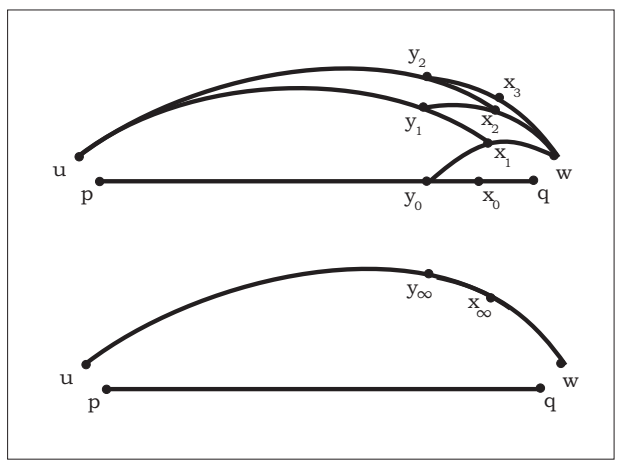}
\end{center}
\caption{Building the sequence of piecewise geodesics, then taking
the limit.} \label{figuw}
\end{figure}

Instead we start with a pair of points $u$ near $\gamma(0)$ and $w$ near $\gamma(T+T_1)$
and define a sequence of piecewise geodesics using the continuous
family about $\gamma([0,T])$ and the short minimal geodesics.  We then prove this carefully
controlled sequence is Cauchy and in fact
converges to a geodesic from $u$ to $w$.  See Figure~\ref{figuw}.

We start with $\sigma_0$ running minimally
from $w$ to $y_0=\gamma(T-T_1)$ and $\gamma_0$ in the given family
from $u$ to the midpoint, $x_1$, of $\sigma_0$.  Then we take
$\sigma_1$
running minimally from $w$ to $y_1=\gamma_0(T-T_1)$ and $\gamma_1$
in the given family from $u$ to the midpoint, $x_2$, of $\sigma_1$,
and so on.
We show both sequences, $\gamma_i$ and $\sigma_i$, are Cauchy using
the
Relative Rauch Comparison Theorem [Theorem~\ref{RelRauch}]
and the fact that their combined
lengths are kept less than $D_\kappa$.  We then prove they converge
to a pair of geodesics $\gamma_\infty$ and $\sigma_\infty$
that glue together to form a geodesic $\gamma_{u,v}$ running from $u$
to $w$.

The proof is very technical because we
control the total height of the bridges as they
build upon one another by choosing very small neighborhoods
for our $u$ and $w$.  We ensure that the accumulated Rauch
estimates and their errors due to the heights remains bounded
by using a geometric series.  

\Pf[\textit{Proof of Theorem~\ref{ultafterpi}}]
Without loss of generality we assume the space is ${\rm CBA}(1)$
and $\gamma$ has length less than $\pi$.  By the paragraph below the
statement
of the theorem, we need only construct a continuous family about
$\gamma$.

Recall the following well-known facts. 
A ${\rm CBA}(1)$ space is locally minimizing
[Proposition~\ref{CBAunifmin}], so it has a positive unique
injectivity radius.  So there exists $T_0>0$ sufficiently small all
geodesics of length less than $ 10 T_0={\rm UniqueInj}(M)$ 
have continuous families about them
[Lemma~\ref{lem-ultinj-fam-2}].

We choose $T_0$ smaller if necesary to ensure
\begin{equation} \label{Tpick1}
\frac{\sin(T_0/6)}{\sin(2T_0/6)} < \frac{3}{4}
\end{equation}
and
\begin{equation} \label{Tpick2}
\cos(T_0/6) + \frac{\sin(T_0)}{\sin(T_0)} < \frac{15}{12}.
\end{equation}
This is possible because the first equation approaches $1/2$ and
the second approaches $7/6$ as $T_0 \to 0$.
Note that for any value of $T\in (T_0, \pi-T_0)$ we then have
\begin{eqnarray} \label{Tpick22}
\frac{|\sin(T-T_0/6)|}{|\sin(T)|}
&=& \frac{|\sin(T) \cos(T_0/6) + \sin(T_0/6) \cos(T)|}{\sin(T)|}\\
&\le& |\cos(T_0/6)| + \frac{|\sin(T_0/6|)}{|\sin(T_0)|}
\frac{|\sin(T_0)|}{|\sin(T)|} |\cos(T)|\\
&\le& |\cos(T_0/6)|
+ \frac{|\sin(T_0/6)|}{|\sin(T_0)|}\frac{|\sin(T_0)|}{|\sin(T)|}
|\cos(T)|\\
&\le &
|\cos(T_0/6)| + \frac{|\sin(T_0/6)|}{|\sin(T_0)|} \cdot 1 \cdot 1 <
\frac{15}{12}.
\end{eqnarray}

Intuitively these equations have been selected to match the
Relative Rauch Comparison (\ref{diff-rel-Rauch})
of Jacobi fields, $J_\sigma$,
running backwards from $\gamma(T+T_0/6)$ so that
\be
\frac{|J_\sigma(T_0/6)|}{|J_\sigma(2T_0/6)|} <
\frac{\sin(T_0/6)}{\sin(2T_0/6)} < \frac{3}{4}
\ee
and Jacobi fields $J_\gamma$ running forwards from $\gamma(0)$,
\be
\frac{|J_\gamma(T-T_0/6)|}{|J_\gamma(T)|}
<\frac{|\sin(T-T_0/6)|}{|\sin(T)|}< \frac{15}{12}.
\ee

Note that $(3/4)(15/12) = 45/48<1$ which will be useful to us
later when we control the accumulation of distances using
a geometric series.

In reality we do not have Jacobi fields and will need to
build bridges and deal with the errors depending on the
height of the bridge given in the Relative Rauch Comparison
Theorem.  Since we will be building bridges upon bridges
in order to construct our geodesic $\gamma$, we need to
control the height and the errors through a series of
estimates.  Nevertheless this choice of $T_0$ suffices
to start our process.

We claim that if we have continuous families about all geodesics
of a given length, $T\in (T_0, \pi-T_0)$,
then we have continuous families about
all geodesics of length $L < T+T_0/6$.  Recall that we do in fact
have
continuous families about geodesics of length $\le 10T_0$ by
our initial choice of $T_0$, so once this claim is proven we can
extend
such continuous families by increments of $T_0/6$ all the way out to
$\pi-5T_0/6$.  So our claim
implies that we can find continuous families for
all geodesics of length less than $\pi-5T_0/6$.  Finally taking $T_0
\to 0$,
we would obtain the result that all geodesics of
length $<\pi$ have continuous families about them.

For now we keep $T_0$ and $T$ fixed and we assume all geodesics of
length
$< T$ have continuous families about them.

Let $\gamma:[0,1]\to X$ have length $L$.

Let $\bar{\gamma}:[0,1]\to X$
be the initial segment of this geodesic running only a length
$T$, so $\bar{\gamma}(s)=\gamma(Ts/L)$.  By the assumption in our
claim,
$\bar{\gamma}$
has a family of geodesics $\bar{F}:U\times V \to \Gamma([0,1],X)$
defined on a neighborhood $U$ of $\bar{\gamma}(0)=\gamma(0)$
and a neighborhood $V$ of $\bar{\gamma}(1)=\gamma(T/L)$
which is continuous at $\bar{\gamma}(0),\bar{\gamma}(1))$.

Let
\begin{equation}
\alpha(h)=\sup_{s\in [T_0/6,\pi-T_0/6]}\alpha_1(s,h)
\end{equation}
where $\alpha_1(s,h)=|s- {\rm arccos}(\cos(s)/\cos(h))|$ as
in the Relative Rauch Comparison Theorem.
Thus $\lim_{h\to 0}\alpha(h)=0$.

In order to create the extended family of geodesics
we need to
insure that we can build bridges between a variety of geodesics
so that we can apply the Relative Rauch Comparison Theorem.  There
exists some $\epsilon_1>0$ such that all bridges of height
$h<\epsilon_1/10$
in a large ball about $\gamma$ will work this way.  More precisely
we choose $\epsilon_1$ small enough that
$h=10\epsilon_1<T_0/4$ can be used as the height
in the Relative Rauch Comparison Theorem so any pair of
geodesics lying within a distance $\epsilon_1$ from ${\gamma}$
form a bridge of height $<h$ with
partition points running through any choice of $r>4h$ and
$R$ as in Theorem~\ref{RelRauch}.
Furthermore we want our bridges to have
controls as strong as (\ref{Tpick1}) and (\ref{Tpick2});
so we require that $h$ be sufficiently small that
$4h<T_0/6$ and
\begin{equation} \label{Tpick1a}
\left( 1+\frac{4h}{D_\kappa} \right)
\frac{\sin(T_0/6)}{\sin(2T_0/6)}
+ \frac{4h+  \alpha(h)}{D_\kappa}
< \frac{3}{4}
\end{equation}
and
\begin{equation}
\left( 1+\frac{4h}{D_{\kappa}} \right)
\left( \cos(T_0/6) + \frac{\sin(T_0/6)}{\sin(T_0)}\right)
+\frac{4h+\alpha(h)}{D_\kappa}
< \frac{15}{12}
\end{equation}
which implies as in (\ref{Tpick22}):
\begin{equation} \label{Tpick2a}
\left( 1+\frac{4h}{D_{\kappa}} \right)
\left(
\frac{|sin(T-T_0/6)|}{|\sin(T)|}
\right)
+\frac{4h+\alpha(h)}{D_\kappa}
< \frac{15}{12}
\end{equation}
for any $T\in (T_0. \pi-T_0)$.

Next we choose  $\delta_1 \in (0,\epsilon_1)$
sufficiently small using the continuity of $\bar{F}$ at
$(\bar{\gamma}(0), \bar{\gamma}(1))$,
such that if $u_1, u_2\in B_{\gamma(0)}(\delta_1)$ and
$v_1,v_2\in B_{\gamma(T/L)}(\delta_1)$
then the geodesic $\bar{F}(u_1,v_1)$ is close to $\bar{F}(u_2,v_2)$:
as a geodesic:
\begin{equation}\label{Tpick4}
d_{\Gamma}(\bar{F}(u_1,v_1),\bar{F}(u_2,v_2))
< (\sum_{j=0}^\infty (45/48)^j)^{-1}(\epsilon_1/4).
\end{equation}
The sum is included here because we will be applying this iteratively
and we need the accumulated error less than $\epsilon_1$ so
that we can apply the Relative Rauch Comparison Theorem to bridges
between geodesics in the family.

Let $\bar{\sigma}:[0,1]\to X$ of length $T_0/3$
be the final segment of the geodesic $\gamma$ running backwards from
$\gamma((T+T_0/6)/L)$ through $\gamma(T/L)$ to
$\gamma((T-T_0/6)/L)$.
Then $\bar{\sigma}$ is minimizing and lies within
a ball of radius $T_0$ about $\gamma(1)$ so it
has a continuous family $F_{min}$ of minimizing
geodesics about it.

There exists a $\delta_2>0$ sufficiently small depending on this
family of minimizing geodesics about $\bar{\sigma}$
such that if
$u_1, u_2\in B_{\sigma(0)}(\delta_2)$ and
$v_1,v_2\in B_{\sigma(1)}(\delta_2)$
then the unique minimizing geodesics $F_{min}(u_i,v_i)$ running
between them
are close together:
\begin{equation}\label{Tpick5}
d_{\Gamma}(F_{min}(u_1,v_1),F_{min}(u_2,v_2))< (\sum_{j=0}^\infty
(45/48)^j)^{-1}(\delta_1/4).
\end{equation}
We further restrict $\delta_2<\epsilon_1/10$ so that there are
bridges between
geodesics which start at a common point in one of the balls
of radius $\delta_2$ and end at points in the other ball of radius
$\delta_2$.

We actually need to refine the neighborhoods two steps further
because the beginning of our iteration does not match the rest of it.

There exists $\delta_3>0$ sufficiently small depending on $\bar{F}$
such that if $u\in B_{\gamma(0)}(\delta_3)$ and
$v\in B_{\gamma(T/L)}(\delta_3)$
then the geodesic $\bar{F}(u,v)$ is close to $\bar{\gamma}$
as a geodesic:
\begin{equation}\label{Tpick6}
d_{\Gamma}(\bar{F}(u,v),\bar{\gamma})
< (\sum_{j=0}^\infty (45/48)^j)^{-1}\delta_2.
\end{equation}

Finally there exists a $\delta_4>0$ sufficiently small depending on
the
family of minimizing geodesics about $\bar{\sigma}$
such that if $u\in B_{\sigma(0)}{(\delta_4)}$ and $v\in
B_{\sigma(1)}(\delta_4)$
then the unique minimizing geodesic $\sigma=F_{min}(u,v)$ running
between
them
is close to $\bar{\sigma}$ as a geodesic:
\begin{equation}\label{Tpick7}
d_{\Gamma}(\sigma,\bar{\sigma})< \delta_3.
\end{equation}

We will now define a map
\begin{equation}
F:B_{\gamma(0)}(\delta_3)\times B_{\gamma(1)}(\delta_4) \to
\Gamma([0,1],X)
\end{equation}
such that $F(\gamma(0),\gamma(1))=\gamma$.
Later we will show it is continuous at $(\gamma(0),\gamma(1))$.
After that we will show it is a continuous family.

Fix $u\in U'=B_{\gamma(0)}(\delta_3)$ and $w \in W'=
B_{\sigma(0)}(\delta_4)=B_{\gamma(T/L)}(\delta_4)$.
We construct a geodesic between them through an iterative limiting
process.

Let
\be
x_0=\gamma(T/L)=\bar{\gamma}(1)=\bar{\sigma}(1/2) \in V
\ee
and let
\be
y_0=\gamma((T-T_0/6)/L)=\bar{\gamma}((T-T_0/6)/T)=\bar{\sigma(1)}.
\ee
Let $\sigma_0$ be the unique minimizing geodesic from $w$ to $y_0$.
It's length is $S_0=d(w,y_0)<T_0$.  Select its midpoint
$x_1=\sigma_0(1/2)$.
By (\ref{Tpick7}) we know
\begin{equation}\label{sigma0}
d_{\Gamma}(\sigma_0,\bar{\sigma})<\delta_3.
\end{equation}
and so $x_1 \in B_{\bar{\gamma}(1)}(\delta_3)$.

Now select a geodesic running from $u$
to $x_1$: $\gamma_0=\bar{F}(u,x_1)$.
Its length is $L_0=L(\gamma_0)=d(u,x_0)$.
We can select a point $y_1=\gamma_0(1-T_0/(6L_0))$ near the endpoint
so that
\begin{equation}
d(y_1,\gamma_0(1))=d(y_0,x_0)=T_0/6.
\end{equation}

By the triangle inequality through $x_1$, we know
$y_1 \in B_{\gamma(1)}(T_0)$.  By (\ref{Tpick6}) and (\ref{sigma0})
we know that in fact
\begin{equation} \label{gamma0}
d_{\Gamma}(\gamma_0,\bar{\gamma})
< (\sum_{j=0}^\infty (45/48)^j)^{-1}\delta_2.
\end{equation}
In particular
\begin{equation} \label{y1y0}
d(y_1,y_0)< (\sum_{j=0}^\infty (45/48)^j)^{-1}\delta_2< \delta_2.
\end{equation}
and
$y_1 \subset
B_{\bar{\sigma}(1)}(\delta_2)=B_{\gamma((T_0-T/6)/L)}(\delta_2)$.

Let $\sigma_1$ be the minimizing geodesic running from $w$ to $y_1$.
Let $x_2$ be it's midpoint.
By (\ref{Tpick5}) and our choice of $\delta_1 \in (0,\epsilon_1)$ we
have
\begin{equation}\label{sigma1}
d_{\Gamma}(\sigma_1,\bar{\sigma})<
(\sum_{j=0}^\infty (45/48)^j)^{-1}(\delta_1/4)<\epsilon_1.
\end{equation}
By our choice of $\epsilon_1$
we can apply the Relative Rauch Comparison
to the bridge between $\sigma_1$ and $\sigma_0$ meeting at
$w$.  By the triangle inequality and the fact that
$d(y_i,\bar{\sigma}(1))<h/2$ and $d(w,\bar{\sigma}(0))<h/2$ we have:
\be
R_{\sigma_i}=d(w,y_i)\in (2T_0/6 - h, 2T_0/6+h)
\ee
and similarly since $d(x_i,\bar{\sigma}(1/2))<h/2$ we have
\be
r_{\sigma_i}\in(T_0/6-h, T_0/6+h).
\ee
So by the Relative Rauch Comparison Theorem and (\ref{Tpick1a}) we
have
\begin{equation}
\frac{d(x_1,x_2)}{d(y_0,y_1)} < \left( 1+\frac{4h}{D_\kappa} \right)
\frac{\sin(T_0/6)}{\sin(2T_0/6)}
+ \frac{4h}{D_\kappa}
\frac{3}{4}.
\end{equation}
In particular
\begin{equation} \label{x2x1}
d(x_2,x_1) < \frac{3}{4}\delta_2
\end{equation}
and by the triangle inequality, (\ref{sigma0}), the choice of
$\delta_3$ and
$\delta_2$ we have
\begin{equation} \label{x2x0}
d(x_2,x_0)< \frac{3}{4}\delta_2+\delta_3
< \frac{\delta_1}{2}+\frac{\delta_1}{2}=\delta_1.
\end{equation}

Now select $\gamma_2=\bar{F}(u,x_1)$.
Its length is $L_1=L(\gamma_1)=d(u,x_1)$
and we can select a point $y_2=\gamma_1(1-T_0/(6L_1))$ near the
endpoint
so that $y_2$ lies near $y_0$.  By (\ref{Tpick4}) and (\ref{x2x0})
we have
\begin{equation}
d_{\Gamma}(\gamma_2,\bar{\gamma})
< (\sum_{j=0}^\infty (45/48)^j)^{-1}(\epsilon_1/4)<\epsilon.
\end{equation}
So we can apply
the Relative Rauch Comparison Theorem to
the bridge formed by $\gamma_1$ and $\gamma_2$ meeting at $u$.
By the triangle inequality and the fact that
$d(u,\bar{\gamma}(0))<h/2$
and $d(x_i,\bar{\gamma}(1))<h/2$ we have
\be
R_{\gamma_i}=d(u, x_{i-1})\in (T-h,T+h)
\ee
and since $d(y_i, \bar{\gamma}((T-T_0/3)/T) )<h/2$ we
have
\be
r_{\gamma_i}=d(u, y_{i-1})\in (T-T_0/3-h,T-T_0/3+h).
\ee
So applying the Relative Rauch Comparison Theorem
and (\ref{Tpick2a})
we have
\begin{equation}
\frac{d(y_2,y_1)}{d(x_2,x_1)}
\left( 1+\frac{4h}{D_{\kappa}} \right)
\left(
\frac{|\sin(T-T_0/6)|}{|\sin(T)|}
\right)
+\frac{4h}{D_\kappa}
< \frac{15}{12}
\end{equation}
In particular
\begin{equation}
d(y_2,y_1) \le \frac{15\cdot 3}{12\cdot 4}  d(y_1,y_0)= \frac{45}{48}
d(y_1,y_0)
\end{equation}
and by (\ref{y1y0}) we have
\begin{equation}
d(y_2,y_0) \le \left(1+ \frac{15}{16}\right) d(y_1,y_0)
< \left(1 + \frac{45}{48} \right) \left(\sum_{j=0}^{\infty} (45/48)^j
\right)^{-1}
\delta_2<\delta_2.
\end{equation}

Continuing iteratively, given $x_0,...x_k$ and $y_0...y_k$
such that each $y_j$ lies on a geodesic $\gamma_{j-1}=\bar{F}(u,x_j)$
running from
$u$ to $x_j$ and each $x_j$ lies on a minimizing geodesic
$\sigma_{j-1}$ running from $w$ to $y_{j-1}$ such that
\begin{equation}
d(x_j,x_{j+1})< \frac{3}{4}d(y_{j-1},y_j)
\end{equation}
and
\begin{equation} \label{79}
d(y_j,y_{j-1}) < \frac{45}{48} d(y_1,y_0).
\end{equation}

We now show we can proceed to defined $x_{k+1}$ and $y_{k+1}$.
Summing the terms in (\ref{79})and applying (\ref{y1y0}) we have
\begin{equation}
d(y_k,y_0) < \sum_{j=0}^k (45/48)^j d(y_1,y_0)<\delta_2.
\end{equation}
So we can draw the next minimizing geodesic $\sigma_{k+1}$
from $w$ to $y_k$, choose its midpoint to be $x_{k+1}$ and
apply the Relative Rauch Comparison
Theorem to get
\begin{equation}
d(x_k,x_{k+1})<\frac{3}{4}d(y_{k-1},y_k)<3\delta_2/4<\delta_1.
\end{equation}
So we can choose $\gamma_{k}=\bar{F}(u,x_{k+1})$, choose
$y_{k+1}=\gamma_{k}((1-T_0/(6L))$, and apply the
Relative Rauch Comparison Theorem to get
\begin{equation}\label{82}
d(y_k,y_{k+1}) <
\frac{15}{12}d(x_k,x_{k+1})<\frac{45}{48} d(y_{k-1},y_k).
\end{equation}
Thus the iteration continues indefinitely.

Notice that the sequence of points $x_k$ is Cauchy
because
\be
d(x_k,x_{k+1})<\left(\frac{3}{4}\right)\left(\frac{15}{12}\right)
d(x_{k-1},x_k).
\ee
So it converges to some point $x_\infty$.
Since $\bar{F}$ is a continuous family on a domain
which contains $x_\infty$, the geodesics
$\gamma_k$ converge to some $\gamma_\infty=\bar{F}(u,x_\infty)$
of length
\begin{equation}\label{lengthinfty}
L_\infty=\lim_{k\to\infty} L(\gamma_k) \in
[L(\bar{\gamma})-\epsilon_1,L(\bar{\gamma})+\epsilon_1]
\end{equation}
because $\gamma_k$ was chosen from the continuous family about
$\bar{\gamma}$
and $\delta_2$ was selected to keep these geodesics much less than
$\epsilon_1$ apart.  In fact
\begin{equation}
d_\Gamma(\gamma_\infty,\bar{\gamma}) <\epsilon_1/2
\end{equation}

The sequence $y_k$ is also Cauchy by (\ref{82})
and converges
to a point $y_\infty$.  Note that
\be
y_\infty=\gamma_\infty (1-T_0/(6L_\infty) ).
\ee
Since minimizing geodesics are continuous with respect to
their endpoints, the $\sigma_k$ also converge to
some minimizing geodesic $\sigma_\infty$ running from
$w$ to $y_\infty$ whose midpoint is $x_\infty$
and length is $T_0/3$.  By our choice of $\delta_3$ we
know that
\begin{equation}
d_\Gamma(\sigma_\infty,\sigma) = \lim_{k\to\infty}
d_\Gamma(\sigma_k,\sigma) < \epsilon_1/2.
\end{equation}

Thus we can build a curve, $\gamma_{u,w}$, running from $u$ to
$w$ first along $\gamma_\infty$ and then backwards along
the second half of $\sigma_\infty$ whose length satisfies
\begin{equation}
L(\gamma_{u,v})=L_\infty+T_0/6\in [L(\gamma)-\epsilon,
L(\gamma)+\epsilon].
\end{equation}
This curve is clearly
locally minimizing about all points in $\gamma_\infty$
and on $\sigma_\infty$ and it is minimzing at the
connection point because it runs minimally from $y_\infty$ to $w$.
Note it is not initially clear that $\sigma_\infty$ and
$\gamma_\infty$ agree between $x_\infty$ and $y_\infty$
but they do have the same length.  Thus
both are minimizing and they must agree since they lie in the
$B_{\gamma(1)}(T_0)$ where minimizing geodesics are unique.
So $\gamma_{u,v}$ is a geodesic running from $u$ to $v$
of length $L(\gamma_{u,v})= L_\infty+T_0/6$.
Furthermore
\begin{equation} \label{dGammauv}
d_\Gamma(\gamma_{u,v}, \gamma) \le
d_\Gamma(\gamma_\infty, \bar{\gamma}) +
d_\Gamma(\sigma_\infty, \bar{\sigma})< \epsilon_1
\end{equation}

Let
\begin{equation}
F:B_{\gamma(0)}(\delta_3)\times B_{\gamma(1)}(\delta_4) \to
\Gamma([0,1],X).
\end{equation}
be defined: $F(u,w)=\gamma_{u,w}$ built as described above.  By its
construction a bridge can be built between $\gamma_{u,w}$
and our original geodesic $\gamma$ restricted to $[0,T+T_0/6]$.

We next need to show that
$F$ is continuous at $(\gamma(0),\gamma(1))$.  Note that the
construction
of $\gamma_{u,v}$ does not depend on the choice of $\epsilon_1$, this
constant was only used to estimate and prove convergence.  If we take
smaller
$\epsilon_1$ we will obtain smaller $\delta_3$ and $\delta_4$, and as
long as we use these smaller values for our domain, we will get
the same values for $F$ just with stronger bounds.

Thus for all $\epsilon>0$, take $\epsilon_1<\epsilon$, determine
$\delta_i$ depending on this value for $\epsilon_1$ and let
$\delta< \min_{i=1,2,3}\{\delta_i\}/\sqrt{2}$.  Then if
\be
(u,v) \in B_\delta(\gamma(0),\gamma(1))\subset X \times X
\ee
we have
\be
u \in B_{\gamma(0)}(\delta_3) \text{ and } v \in
B_{\gamma(0)}(\delta_2)
\ee
so by (\ref{dGammauv})
\be
d_\Gamma(F(u,v),\gamma)< \epsilon_1 <\epsilon.
\ee

Thus we have shown that if we have families about
all geodesics of length $L < T+T_0/6$ which are continuous.

Thus we have shown that all geodesics $\gamma$ of length
$L <T +T_0/6$ have families defined about
them which are continuous at $(\gamma(0),\gamma(1))$.
By Lemma~\ref{lem-ult4}, we thus have our claim that
all geodesics $\gamma$ of length
$L <T +T_0/6$ have continuous families defined about them.
The theorem follows.
\qed

\section{Applications to ${\rm CBA}(\kappa)$ spaces}\label{sectappl}

In this section we briefly survey some applications
of our results on geodesic spaces
to ${\rm CBA}(\kappa)$ spaces and potential further
directions of research.  We begin with the long homotopy
lemma and then the injectivity radius theorem
and finally discuss some further directions.

\begin{theo} \label{longhomCBA}
If $M$ is a locally compact ${\rm CBA}(\kappa)$ space and $c:[0,1]\to M$
is a nontrivial contractible closed geodesic of
length $l(c)< 2D_\kappa$, then any null homotopy for
$c$ contains a curve of length $\ge 2 D_\kappa$.
\end{theo}

\Pf
We begin by noting that $M$ is locally uniformly minimizing
[Lemma~\ref{CBAunifmin}].
By Theorem~\ref{ultafterpi} we know ${\rm conj} (M)\ge D_\kappa$.
Thus this follows immediately from Theorem~\ref{longhomotopy}, our
generalization of the long homotopy lemma.
\qed

The standard application of the theorem to Riemannian manifolds is
depicted in Figure~\ref{longhom}.  An application to a ${\rm CBA}(\kappa)$
space can be seen by taking two copies of that example and joining
them together at finitely many corresponding points similar to
the Tetra bi-sphere (Example~\ref{tetrasphere}).  Note how it is
crucial in this application that the fans constructed do not
require uniqueness of geodesics between the points.

\begin{problem} \label{longhom-morse}
The Riemannian long homotopy lemma has also often been used in
combination with Morse Theory to prove the existence of smooth
closed geodesics.
Possible extensions of Morse theory to length spaces appear in
\cite{SorLength}.
Thus one might try to extend some of this existence theory
to length spaces and ${\rm CBA}(\kappa)$ spaces.
\end{problem}

Our next application is closely related to a step in Charney-Davis'
proof of Gromov's Systole Theorem \cite{chardavis} (c.f.\ \cite{BH}, II.4.16).

\begin{coro}[Charney--Davis] \label{WASCOR10.8}
A compact ${\rm CBA}(\kappa)$ space, $X$, with $\kappa>0$
such that ${\rm MinRad}(X)\in (0,D_\kappa)$ has a closed
geodesic with length twice the
minimal injectivity radius.
\end{coro}

\Pf
We just combine Theorem~\ref{kling2} with Theorem~\ref{ultafterpi}.
\qed

In the Charney-Davis proof of Gromov's systole theorem
they essentially prove the following statement:
{\em If the Unique Injectivity radius is less than $D_\kappa$ then there
is a closed geodesic of length equal to twice the unique injectivity
radius which is in fact a digon, so its length is twice the minimal
injectivity radius as well.}  They do not explicitly state this but the essence
of the idea is there so we consider Corollary~\ref{WASCOR10.8} to
already be known.  The argument in their situation is easier than ours
because it is easy to construct continuous families about unique geodesics.
So Corollary~\ref{WASCOR10.8} is an over simplification of our Klingenberg
Injectivity Radius Theorem (Theorem~\ref{kling2}) but our theorem
holds in a much wider setting.

We close with a discussion of the implications of equality in
Theorem~\ref{ultafterpi}.

\begin{defn} \label{defnrank1}
A ${\rm CBA}(1)$ space is said to have positive spherical rank
if any geodesic segment of length $\pi$ has conjugate endpoints.
\end{defn}

On Riemannian manifolds, one has the immediate consequence
that every geodesic has a Jacobi field running along it
of the form $J(t) = \sin(t)E(t)$ where $E$ is parallel.  In
\cite{ssw}, the authors called
closed manifolds with this property to have
\textit{positive spherical rank}. In that paper it is shown:
{\em If a Riemannian manifold
of sectional curvature at most 1 has positive spherical rank,
then its universal cover is isometric to a compact, rank one,
symmetric space}.

\begin{problem} \label{problem-rank1}
Can one classify ${\rm CBA}(\kappa)$ spaces with positive spherical rank?
\end{problem}

Note that following example demonstrates that the classification
cannot be restricted to symmetric spaces at least without the further
assumption of a lower curvature bound.

\begin{ex}\label{triplehemisphere}
The triple hemisphere is a compact length space created by
gluing three hemispheres together along a common equator.
It is a ${\rm CBA}(1)$ space with infinite negative curvature
along the equator.  See for example \cite{BBI} for theorems
about constructing ${\rm CBA}(\kappa)$ spaces using gluing.
One can find the distance between any pair of points on
a common hemisphere using standard spherical geodemtry and
the distance between any pair of points on different
hemispheres by considering the two hemispheres as a single sphere.

Any geodesic of length $<\pi$ will pass through at most two
hemispheres and so can be seen to be running along a great
circle on the sphere formed by those two hemispheres.  Thus
it has no ultimate conjugate points before $\pi$ (otherwise the
sphere would have a conjugate point by Proposition~\ref{prop-ult4}.
So the triple hemisphere is a space with positive spherical rank.
\end{ex}

While the above example is not a symmetric space, it is an example of a
\textit{spherical building}. One might hope that a classification of spaces with
positive spherical rank might be limited to spherical buildings (or their generalizations)
but A.\ Lytchak informs us that in dimensions 3 and higher it is possible to
construct fairly complicated spaces that have positive spherical rank. In the 2-dimensional
case, W.\ Ballmann and M.\ Brin showed that if one further assumes that the space is
obtained by gluing spherical simplices together then one has rigidity. Specifically they
showed that all such examples are spherical buildings, spherical joins or one class
of examples that is neither; see \cite{bb}. One can see that the examples in \cite{bb} that are neither
buildings nor joins do not have positive spherical rank. So in this case one may be able
to prove that such spaces are rigid.

\appendix
\section{Translated Excerpt from Rinow \cite{Rinow}}

{\em As Rinow's text is in German, we have translated relevant
sections here and added comments in italics}.

\vspace{.2cm}

\subsection{Absolute Conjugate Points; page 172}

\vspace{.2cm}

\noindent
{\em Rinow is first interested in minimizing geodesics}:

\vspace{.2cm}

Let $f(s), 0 \leq s < \beta$ denote the normal representation 
of a geodesic ray $S_a$ with starting point $a = f(0)$. Let 
$k(S_a)$ denote the supremum over all $s' \in [0, \beta)$ for those 
$s'$ for which $f\mid_{[0,s']}$ is a shortest curve (i.e., geodesic 
segment) between $a$ and $f(s')$. Obviously, $0 < k(S_a) \leq \beta$. 

\vspace{.2cm}

\noindent
{\em Rinow then defines a notion of conjugate point which is
unrelated to our notion and not equivalent to the Riemannian
notion}:

\vspace{.2cm}

The {\bf absolute conjugate point of the ray $S_a$} is defined to be 
the point $f(k(S_a))$ if $k(S_a) < \beta$, or in the case 
$k(S_a) = \beta < \infty$ the limit point of $S_a$ if it exists. 
We also denote it as $f(k(S_a))$. In case $k(S_a) < \beta$, $f(s)$ 
always represents a shortest curve on $[0, k(S_a)]$, however, not 
on any interval $[0,\alpha]$, where $k(S_a) < \alpha \leq \beta$. 
It follows that all inner points for $f(k(S_a))$ are shortest 
curves for $S_a$. In an inner metric space, obviously 
$k(S_a) = d(a, f(k(S_a)))$, if $f(k(S_a))$ always exists.

\vspace{.2cm}

\noindent
{\em Rinow later defines a notion equivalent to a one sided
conjugate point:}

\vspace{.2cm}

\subsection{Conjugate Points pp 414-415}


\vspace{.2cm}

\noindent $X$ is an inner metric ({\em geodesic})
space with the following properties (page 414):
\smallskip

\noindent (a) $X$ is locally compact.
\smallskip

\noindent (b) $X$ does not have any branching points (\textit{i.e., $X$ does not have bifurcation of geodesics; see page 162}).
\smallskip

\noindent (c) For each point $a \in X$, there exists $\epsilon_a >0$ so that
$$
\text{inf } \{ k(y): y \in U(a,\epsilon_x)\} > 0.
$$

\vspace{.2cm}

\noindent
{\em Note that due to the nonbranching, there is local uniqueness
of geodesics of length less than $k(y)$.  One can also piece together
geodesics to extend them:}

\vspace{.2cm}

Let $a\in X$ be any point. Let $\alpha$ be a real number 
so that $0 < \alpha < k(a)$. Let $\sum_\alpha$ denote the 
peripheral sphere of $a$ with radius $\alpha$ (points in $X$ at 
distance $\alpha$ to $a$). Every point $\xi \in \sum_\alpha$ can 
be connected to $a$ by exactly one shortest curve and there 
is exactly one geodesic ray $S_\xi$ which is infinitely 
extendable. 

\vspace{.2cm}

\noindent
{\em Rinow then defines an exponential map, $f$, using this uniqueness:}

\vspace{.2cm}

Let
$f(\xi, s), 0\leq s < \infty$ be the normal parameter representation of $S_\xi$. We set $P = \sum_\alpha \times [0, \infty)$ with the product metric
$$
\rho((\xi, s), (\xi', s')) = \sqrt{d(\xi, \xi')^2 + |s-s'|^2}
$$
Then $f: P \rightarrow X$ is a well defined map. Every point $x\neq a$ can be connected to $a$ be a shortest curve. Consequently, we get a geodesic ray through $x$. Moreover, $f(\xi,0) = a$ for every $\xi \in \sum$ and $f$ is continuous on $P$.

\vspace{.2cm}

\noindent
{\em Rinow uses this notion to define conjugate points which
he next proves is equivalent to what we call a one sided conjugate point:}

\vspace{.2cm}

Given a point $(\xi_0, s_0) \in P$ with $s_0 >0$, we say $s_0$ is 
an {\bf ordinary point} for the ray $f(\xi_0, s_0)$ if for 
some $\delta >0$, $f$ restricted to 
$V_\delta (\xi_0, s_0) = \{(\xi, s): 
d(\xi, \xi') < \delta, |s-s_0| < \delta \}$, is a topological map 
(i.e., homeomorphism onto its image see p 11) and 
in addition $f(V_\delta (\xi_0, s_0))$ is open in $X$. If no 
such $\delta >0$ exists, then $s_0$ is called a 
{\bf conjugate point}. Obviously, the set of ordinary points 
is open while the set of conjugate points is closed in $[0, \infty)$.

\vspace{.2cm}

\noindent
{\em Rinow then proves the following theorem which we only state here:}

\vspace{.2cm}

\begin{theorem} {\em 
Let $X$ be an $n$-dimensional manifold. $s_0 >0$ is an ordinary point 
of $f(\xi_0, s)$ if and only if there exists $\delta >0$ such 
that $f$ restricted to $V_\delta (\xi_0, s_0)$ is injective. }
\end{theorem}

\vspace{.2cm}

\noindent
{\em One then sees that this is equivalent to the notion of a one-sided 
conjugate point because $f$ was an exponential map.  That is, if
$s_0$ is an ordinary point of the exponential map $f$, then
$\xi_0(s_0)$ is not a one sided conjugate point of $p$, otherwise
the existence of $\gamma_i$ and $\sigma_i$ would contradict the
local injectivity of $f$ for $i$ sufficiently large and visa versa.}

\medskip

%

\bibliographystyle{alpha}

Department of Mathematics, University of Oklahoma, shankar@math.ou.edu 

Lehman College and CUNY Graduate Center, sormanic@member.ams.org

\end{document}